\documentclass[12pt,reqno,a4paper]{gen-j-l}
\usepackage{a4wide,epsfig,amssymb,amsmath}

\newtheorem{thm}{Theorem}[section]
\newtheorem{lem}[thm]{Lemma}

\newtheorem{coro}[thm]{Corollary}
\newtheorem{prop}[thm]{Proposition}

\newcommand{\inte}{\operatorname{int}}
\newcommand{\Bd}{\partial}

\newcommand{\Z}{\mathbf{Z}}
\newcommand{\N}{\mathbf{N}}

\begin{document}
\title{Parameterizations of 1-bridge torus knots}
\author{Doo Ho Choi}
\email{dhchoi@knot.kaist.ac.kr}
\address{Department of Mathematics, Korea Advanced Institute
of Science and Technology, Taejon, 305-701, Korea}
\author{Ki Hyoung Ko}
\email{knot@knot.kaist.ac.kr}
\address{Department of Mathematics, Korea Advanced Institute
of Science and Technology, Taejon, 305-701, Korea}
\begin{abstract}
A 1-bridge torus knot in a 3-manifold of genus $\le 1$ is a knot
drawn on a Heegaard torus with one bridge. We give two types of
normal forms to parameterize the family of 1-bridge torus knots
that are similar to the Schubert's normal form and the Conway's
normal form for 2-bridge knots. For a given Schubert's normal form
we give algorithms to determine the number of components and to
compute the fundamental group of the complement when the normal
form determines a knot. We also give a description of the double
branched cover of an ambient 3-manifold branched along a 1-bridge
torus knot by using its Conway's normal form and obtain an
explicit formula for the first homology of the double cover.

\end{abstract}

\maketitle

\section{Introduction}
One of traditions in knot theory is to study a family of knots
satisfying a certain condition. Examples of such families include
the family of torus knots studied by Dehn and Schreier and the
family of 2-bridge knots studied by Schubert, Montesinos and
Conway. These classes can be referred as the classes of knots and
links indexed by the pairs $(g,b)$ of non-negative integers as
defined in \cite{doll}. A knot $K$ in a 3-manifold $M$ has a {\em
$(g,b)$-decomposition} or is called a {\em $(g,b)$-knot} if for
some Heegaard splitting $M=U\cup V$ of genus $g$, each of $K\cap
U$ and $K\cap V$ is consisted of trivial $b$ arcs. A collection of
properly embedded arcs in a 3-manifold $W$ with boundary is {\em
trivial} if arcs in the collection together with arcs on
$\partial W$ joining the two ends of the arcs bound mutually
disjoint disks, called {\em cancelling disks}, in $W$. A
$(g,b)$-knot can be embedded in a Heegaard surface of genus $g$ in
$M$ except at $b$ over(or under)-bridges and vice versa. Torus
knots are $(1,0)$-knots and 2-bridge knots are $(0,2)$-knots.
Clearly the family of $(g,b)$-knots becomes strictly larger as $g$
or $b$ increases. Since an over-bridge can be removed by adding a
handle and by embedding the over-bridge into the added handle,
$(g,b)$-knots are contained in the family of $(g+1,b-1)$-knots.

In this paper we study 1-bridge torus knots, that is,
$(1,1)$-knots in a 3-manifold. An ambient 3-manifold necessarily
has the Heegaard genus $\le 1$ and so it can be either $S^3$ or a
lens space. The family of 1-bridge torus knots in $S^3$ contains
torus knots and 2-bridge knots and is contained in the family of
double torus knots, that is, $(2,0)$-knots. Hill and Murasugi
studied the family of double torus knots in \cite{hill,murasugi}
and parametrized the family. Non-trivial knots with the trivial
Alexander polynomial was found in the subfamily of double torus
knots that separate the double torus. They also considered
non-separating double torus knots and a subfamily of 1-bridge
torus knots and found various double torus knots that are fibered.

Every 1-bridge torus knot has the tunnel number one, but not all
tunnel-number-one knots are 1-bridge torus knots. In
\cite{mori-saku-yoko1}, Morimoto, Sakuma and Yokota found
tunnel-number-one knots that are not 1-bridge torus knots as
confirmed by a condition on the Jones polynomial for a knot to
admit a $(g,b)$-decomposition in \cite{yokota}. In
\cite{mori-saku-yoko2}, they gave another criteria to determine
whether a given knot has the tunnel number one and whether it is a
1-bridge torus knot.

Besides torus knots and 2-bridge knots, the family of 1-bridge
torus knots includes Berge's double-primitive knots, 1-bridge
braids that classified by Gabai in \cite{gabai}, and satellite
knots of tunnel number one. Morimoto and Sakuma studied satellite
knots of tunnel number one and classified their unknotting tunnels in
\cite{mori-saku}.

Fujii showed that any Alexander polynomial can be realized by a
1-bridge torus knot in \cite{fujii}. He also found a family of
non-trivial 1-bridge torus knots with trivial Alexander
polynomial. In \cite{choi-ko}, a subfamily of 1-bridge torus knots
is completely classified using their genera and Jones polynomial.
Also cyclic branched covers of ambient spaces along 1-bridge torus
knots are known as Dunwoody 3-manifolds and are studied in
\cite{song-kim} and \cite{gra-mul}.

\begin{figure}[!ht]
$$\epsfbox{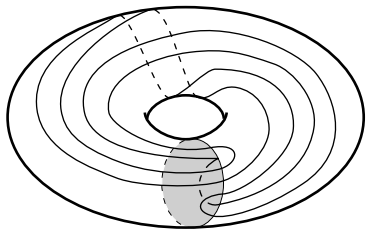}$$
\caption{1-bridge torus knot}\label{fig:exam1}
\end{figure}

We will parameterize the family of 1-bridge torus knots using two
kinds of normal forms as done for the family of 2-bridge knots.
Schubert described a 2-bridge knots by a pair of integers of a
certain condition from its top view. In the top view a 2-bridge
knots is embedded in a plane except the two bridges. He in fact
completely classified 2-bridge knots using this normal form
\cite{schubert}. Since a 1-bridge torus knot can be embedded in a
standard torus except the bridge (See Figure~\ref{fig:exam1}), we
will describe it by a 4-tuple of integers from this top view.
Section 2 will be devoted to simplify an embedded curve with two
fixed ends on a torus up to surface isotopies that disturb neither
bridges nor knot types of 1-bridge torus knots. We show that a
4-tuple of integers is enough to describe such a curve and give an
algorithm to count the number of component of the curve given by a
4-tuple. We will call such a 4-tuple the \emph{Schubert's normal
form} of the 1-bridge torus knot determined by a 4-tuple. The
Schubert's normal form is useful to compute the fundamental group
of the exterior of a 1-bridge torus knot. In Section 3, we give an
algorithm to compute the knot group from the Schubert's normal
form. As a corollary, we determine when the exterior of a 1-bridge
torus knot in a lens space has a double cover.

On the other hand, a 2-bridge knot can also be viewed as a 4-plats
as studied first in \cite{schumann}. This corresponds to a side
view and the composition of homeomorphisms of a four-punctured
sphere determines the 2-bridge knot. Using this description,
Conway constructed a bijection between 2-bridge knots and lens
spaces via double branched covers \cite{conway}. A similar
description using the composition of homeomorphisms on a
two-punctured torus is possible for 1-bridge torus knots and this
will be called the \emph{Conway's normal form}. In Section 4, we
describe a Conway's normal form of a 1-bridge torus knot and we
show it is well-defined in the sense that a Conway's normal form
belongs to a free subgroup generated by three homeomorphisms in
the mapping class group of a two-punctured torus. Finally we
construct the double branched cover of an ambient space branched
along a 1-bridge torus knot given by the Conway's normal form and
give a formula for the first homology of the double branched
cover.

\section{An embedded curve in a torus}
\begin{figure}[!ht]
$$\epsfbox{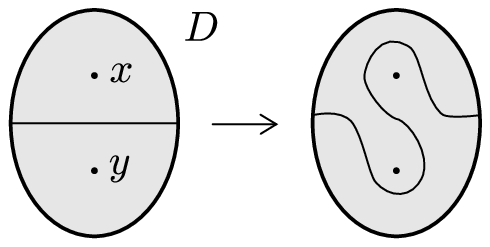}$$
\caption{}\label{fig:exchange}
\end{figure}

Since 1-bridge torus knots can be embedded in a standard torus $T$
except the bridge connecting two points $x, y$ in $T$, we will try
to classify embedded curves with two fixed ends $x, y$ in $T$ up
to isotopies of two-punctured torus that preserves the bridge and
the knot type of a given 1-bridge torus knot and extends to
ambient isotopies on $S^3$.
Let $m$ be a meridian curve on $T$ containing the points $x, y$.
Consider the following two types of isotopies on $T$.
\begin{enumerate}
\item[(I)] An isotopy $h_t$ fixing $x$, $y$ pointwise such that $h_0$ is an
identity map;
\item[(II)] An isotopy $h_t$ such that $h_0$ is an identity map and
$h_1$ is a homeomorphism exchanging $x$ and $y$ counterclockwise
or clockwise and fixed on $T-D$ for a small disk $D$ containing
$x, y$ as illustrated in Figure~\ref{fig:exchange}.
\end{enumerate}

\subsection{A normal form of an embedded curve}\label{subsec:simple}
\begin{figure}[!ht]
$$\epsfbox{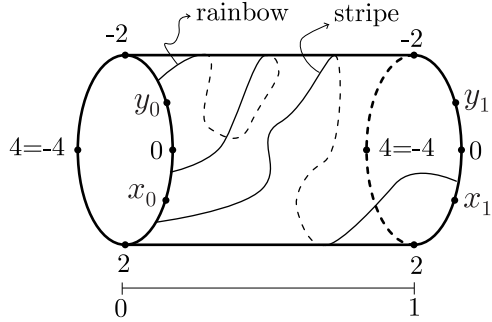}$$
\caption{Two types of arcs on the cylinder}\label{fig:arc-type}
\end{figure}

We parameterize the torus $T$ as $S^1\times[0,1]/x\times 1\sim
x\times 0$ for $x\in S^1=[-4,4]/4\sim-4$. Let $x,y$ be $[1\times
0], [-1\times 0]$, respectively. Let $\gamma$ be an embedded curve
on $T$ with two ends $x, y$. We suppose that $\gamma$ is not
isotopic to an arc on $S^1\times 0\subset T$ and $|m\cap \gamma|$
is minimal up to isotopies of type I and II. If we cut the torus
$T$ along the meridian circle $m=S^1\times 0$, we obtain a
collection $\Gamma$ of arcs from $\gamma$ on the cylinder
$C=S^1\times I$. We may assume that $\Gamma$ contains no closed
curves since we are not interested in trivial components that
splits. Then each arc $\alpha\in \Gamma$ is one of the following
two types (See Figure~\ref{fig:arc-type}):

\begin{enumerate}
\item An arc $\alpha$ is called a \emph{rainbow} if either $\Bd\alpha\cap(S^1\times0)=\emptyset$ or
$\Bd\alpha\cap(S^1\times1)=\emptyset$.
\item An arc $\alpha$ is called a \emph{stripe} otherwise.
\end{enumerate}
Let $x_0=(1,0)$, $y_0=(-1,0)$, $x_1=(1,1)$, and $y_1=(-1,1)$ in
$C$. If $\alpha$ is a rainbow in $C$, then $C-\alpha$ has two
components, one is a disk and the other is a cylinder. We denote
the disk by $D_\alpha$. If $D_\alpha$ contains a point $t$ in
$S^1\times 0$ or $S^1\times 1$, we say that \emph{$\alpha$
contains $t$}. Important properties of the set $\Gamma$ of arcs
are collected in the following lemma.

\begin{lem}\label{lem:xy-stripe} Under the above assumtion on
$\gamma$ and $\Gamma$, the following holds:
\begin{enumerate}
\item A rainbow contains one and only one of $x_0,
y_0, x_1, y_1$.
\item There are no two arcs starting from both $x_0$ and
$y_0$ nor both $x_1$ and $y_1$.
\item An arc starting from one of $x_0, x_1, y_0, y_1$
is a stripe.
\item If there is a stripe starting from $x_0$(resp. $x_1$) then
there is also a stripe starting from $y_1$(resp. $y_0$).
\item Each end of a rainbow is joined to an end of a stripe.
\end{enumerate}
\end{lem}
\begin{proof} (1) If there is a rainbow containing either both
$x_0$ and $y_0$ or both $x_0$ and $y_0$, $\gamma$ joining $x$ and
$y$ can be isotoped into the meridian $m$ and this violates our
assumption. If there is a rainbow containing none of $x_0, x_1,
y_0, y_1$, then the rainbow can be remove from $\Gamma$ by
isotopies of type I and this contradicts the minimality of $|m\cap
\gamma|$.

(2) Suppose that $\Gamma$ has two stripes starting from both $x_0$
and $y_0$, respectively. Since $\gamma$ joins $x$ and $y$ on $T$
the two stripes are connected via since $\gamma$ joins $x$ and $y$
and then $\Gamma$ must have a rainbow whose ends belong to
$S^1\times 0$. Since the rainbow contains none of $x_0$ and $y_0$,
$|m\cap l|$ can be reduced.

(3)  Let $\alpha$ be an arc in $\Gamma$ starting from $x_0$.
Suppose $\alpha$ is a rainbow. Then $\alpha$ must contain $y_0$.
Otherwise $|m\cap \gamma|$ can be reduced by an isotopy of type
(I).
\begin{figure}[!ht]
$$\epsfbox{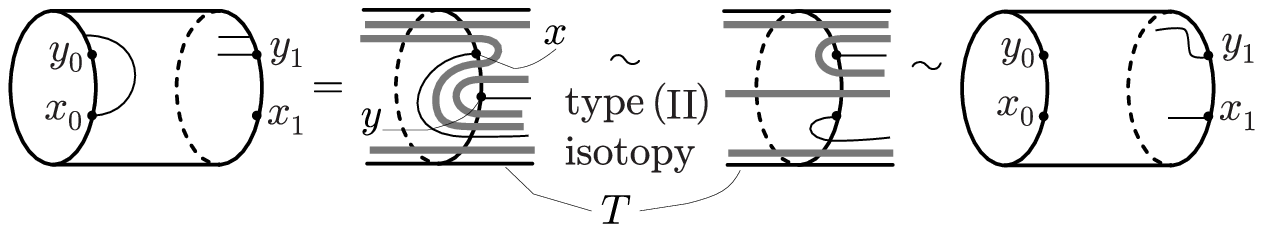}$$
\caption{}\label{fig:iso-disk}
\end{figure}
By an isotopy of type II as illustrated in
Figure~\ref{fig:iso-disk}, $\alpha$ can be removed from $\Gamma$
and $|m\cap \gamma|$ was not minimal.

(4) Since $x$ and $y$ are connected via $\gamma$, this immediately
follows from (2) and (3).

(5) Suppose that there is a rainbow whose one end is connected
with that of the other rainbow. Then we have the two cases as in
Figure~\ref{fig:2-cases-r}.
\begin{figure}[!ht]
$$\epsfbox{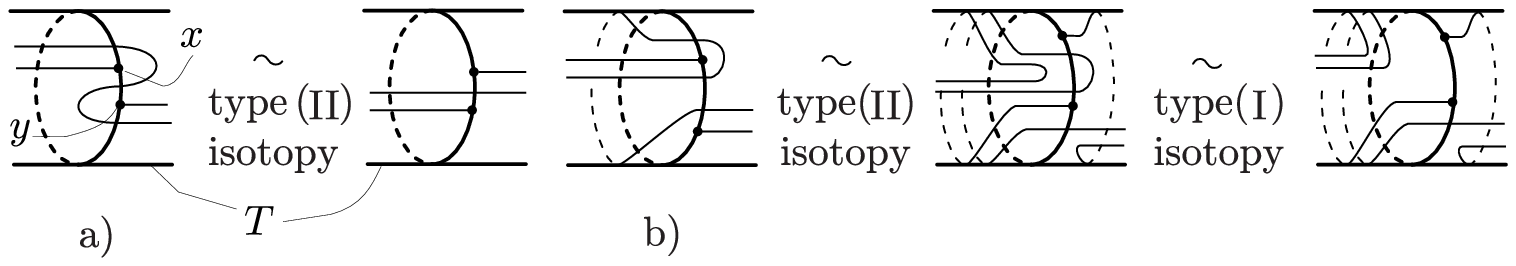}$$
\caption{}\label{fig:2-cases-r}
\end{figure}
For each case of Figure~\ref{fig:2-cases-r}, we can reduce $|m\cap
\gamma|$ by isotopies of type I and II as shown in
Figure~\ref{fig:2-cases-r}.
\end{proof}

If $L$ has a stripe starting from $x_0$(resp. $y_0$), $\gamma$ is
called \emph{$+1$-type}(resp. \emph{$-1$-type}).

\begin{lem}\label{lem:n-rainbow}
Any rainbow contains $x_1$ or $y_0$(resp. $x_0$ or $y_1$) if
$\gamma$ is $+1$-type(resp. $-1$-type). And the number of rainbows
containing $x_0$(resp. $y_0$) is equal to that of rainbows
containing $y_1$(resp. $x_1$).
\end{lem}
\begin{proof}
Suppose $\gamma$ is $+1$-type. Then a rainbow can not contain
$x_0$ nor $y_1$ and if it contains none of $x_1$, $y_0$ then it
can be removed by an isotopy of type I but it is impossible since
$|m\cap \gamma|$ is minimal. Therefore it must contain $x_1$ or
$y_0$. And since $\gamma$ is embedded, the number of rainbows
containing $x_1$ is equal to that of rainbows containing $y_0$.
\end{proof}

\begin{thm}\label{thm:arc-class}
Let $\gamma$ be an embedded curve on the torus $T$ with $\partial
\gamma=\{x, y\}$ and $|m\cap \gamma|$ be minimal up to isotopies
of type I and II where $m$ is a meridian circle containing $x, y$.
Then $\gamma$ is represented by a 4-tuple
$(r,s,t,\rho)_{\epsilon}$, where $r,s,t$ are non-negative
integers, $\rho$ is an integer and $\epsilon$ is a sign $\pm 1$.
If $s\neq 0\neq t$ then $(r,s,t,\rho)_{\epsilon}$ is unique and if
$s$(resp. $t$) is $0$ then a 4-tuple $(r,0,t,\rho)_{-\epsilon}$
(resp. $(r,s,0,\epsilon2(2r+1)+\rho)_{-\epsilon}$) also represents
$\gamma$.
\end{thm}
\begin{proof}
\begin{figure}[!ht]
$$\epsfbox{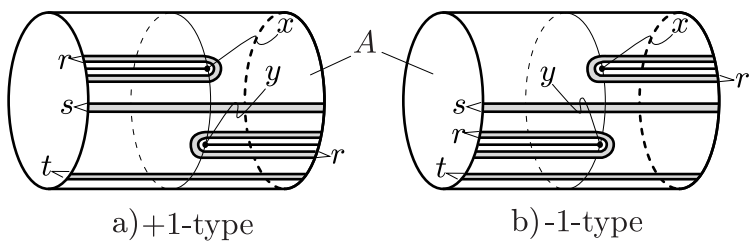}$$
\caption{Shape of $\gamma$ on $A$}\label{fig:shape-arc}
\end{figure}
Let $\gamma$ be an embedded curve of $\epsilon$-type. From the
view of Lemma~\ref{lem:xy-stripe} and Lemma~\ref{lem:n-rainbow},
$\gamma$ must look like Figure~\ref{fig:shape-arc} on a
neighborhood $A$ of the meridian $m$.

If we identify two ends of $A$ by the suitable
$\frac{2\pi\rho}{n}$-rotation map on $S^1$ then we have an
embedded curve isotopic to $\gamma$ on $T$ where $n=2r+1+s+t$.
Therefore $(r,s,t,\rho)_\epsilon$ represents $\gamma$ where
$r,s,t$ are nonnegative integers, $\rho$ is an integer and
$\epsilon=\pm 1$.

The 4-tuple $(r,s,t,\rho)_{\epsilon}$ is unique up to isotopies of
type I and II if $s\neq 0\neq t$ and $|m\cap \gamma|$ is minimal.
And if $s$ or $t$ is $0$ then the last statement of the theorem
holds by exchanging rainbows via isotopies of type II.
\end{proof}

\subsection{Components counting algorithm}\label{subsec:comp}
Given non-negative integers $r,s,t$, integer $\rho$ and a sign
$\epsilon$, we will try to find a condition that makes a 4-tuple
$(r,s,t,\rho)_\epsilon$ represent a simple arc on $T$ with ends
$x,y$. We will in fact give an algorithm to count the number of
components of the curve determined by a 4-tuple.
\begin{figure}[!ht]
$$\epsfbox{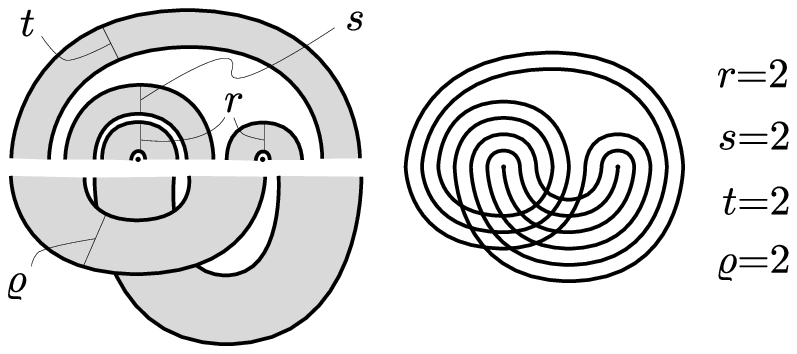}$$
\caption{}\label{fig:comp1}
\end{figure}
We consider a set of arcs on a cylinder such as
Figure~\ref{fig:shape-arc} determined by $r,s,t$ and $\epsilon$.
Let $n=2r+1+s+t$. If we identify two ends of the cylinder by a
$\frac{2\pi\rho}{n}$ rotation, then we obtain an embedded curve on
a torus $T$. Denote the number of components of this curve by
$|(r,s,t,\rho)_\epsilon|$. It is easy to see that
$|(r,s,t,\rho)_{+1}|$ is equal to the number of components of the
curve depicted in Figure~\ref{fig:comp1}. We first remark that
\begin{enumerate}
\item $|(r,s,t,\rho)_{-1}|=|(r,s,t,n-\rho)_{+1}|$.
\item $|(r,s,t,\rho)_\epsilon|=|(r,s,t,\bar{\rho})_\epsilon|$,
where $\bar{\rho}\equiv \rho\mod n$.
\end{enumerate}
Thus we assume that $\epsilon=+1$ and $0\leq \rho<n=2r+1+s+t$.
\begin{figure}[!ht]
$$\epsfbox{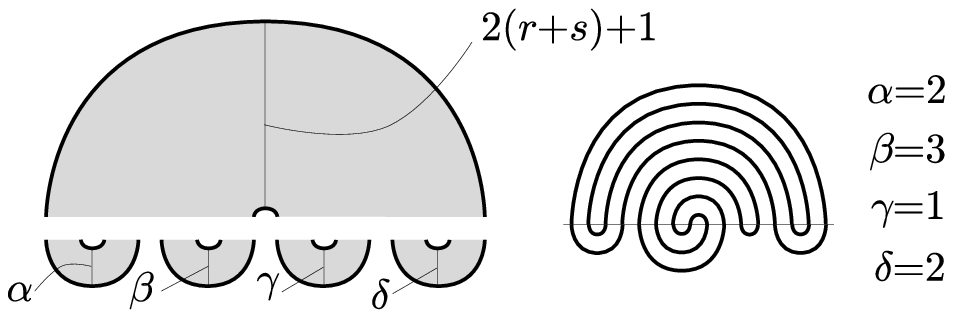}$$
\caption{}\label{fig:comp2}
\end{figure}
\begin{prop}\label{prop:cc-method}
$|(r,s,t,\rho)_{+1}|$ is equal to the number of components in the
curve depicted in Figure~\ref{fig:comp2}, where if $0\le \rho <
n-t$ then $\beta\equiv t \mod \rho$ with $0\le\beta<\rho$,
$\alpha=\rho-\beta$, $\gamma=n-\rho-t$, and $\delta=s$, or if
$n-t\le \rho < n$ then $\beta\equiv t \mod (n-\rho)$ with
$0\le\beta<n-\rho$, $\alpha=\rho-t$, $\gamma=n-\rho-\beta$, and
$\delta=s$.
\end{prop}
\begin{proof} In the figures of the proof, two curves being
``$=$'' mean that the curves have the same number of components.\\
(i) If $t<\rho$ and $t<n-\rho$, we may isotope the upper rainbow
with $t$ strands as in Figure~\ref{fig:equi-1} to remove it.
\begin{figure}[!ht]
$$\epsfbox{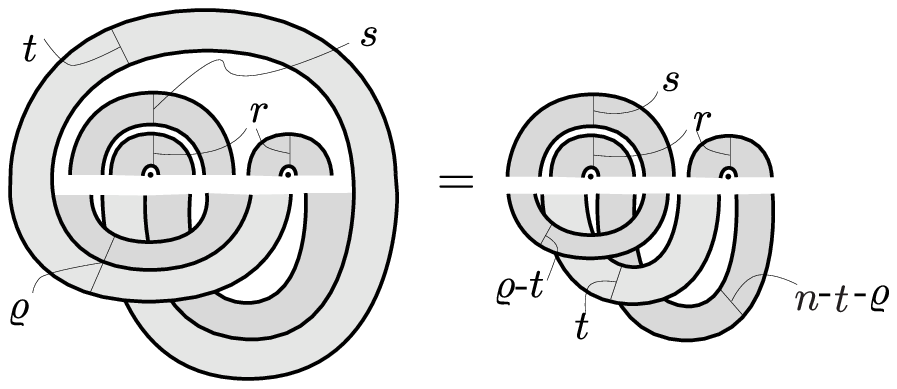}$$
\caption{}\label{fig:equi-1}
\end{figure}

(ii) If $t\geq\rho$, we may isotope parts the upper rainbow
consisted of $\rho$ strands at a time until it becomes the case
(i) as in Figure~\ref{fig:equi-2} where $\beta$ is the remainder
when $t$ is divided by $\rho$.
\begin{figure}[!ht]
$$\epsfbox{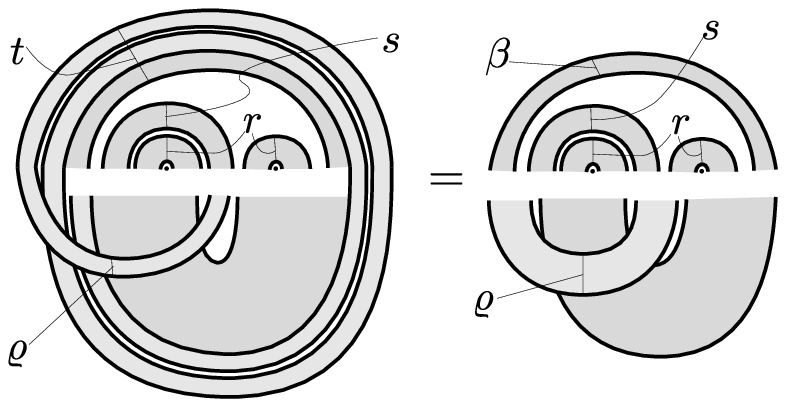}$$
\caption{}\label{fig:equi-2}
\end{figure}

(iii) If $t\geq n-\rho$, we similarly isotope parts of the upper
rainbow consisted of $n-\rho$ strands at a time until it becomes
the case (i) where $\beta$ is the remainder in the division of $t$
by $n-\rho$.

In the all three cases $|(r,s,t,\rho)_{+1}|$ is equal to the
number of components of the curve on the left in
Figure~\ref{fig:equi-4} where if $0\le \rho < n-t$ then
$\beta\equiv t \mod \rho$ with $0\le\beta<\rho$,
$\alpha=\rho-\beta$, and $\gamma=n-\rho-t$, or if $n-t\le \rho <
n$ then $\beta\equiv t \mod (n-\rho)$ with $0\le\beta<n-\rho$,
$\alpha=\rho-t$, and $\gamma=n-\rho-\beta$. In
Figure~\ref{fig:equi-4}, the curve in the middle is obtained from
the one on the left by creating a large curl involving $s$ strands
that maintains the connections on strands involved. The curve on
the right is obtained from the one in the middle by unfolding.
\begin{figure}[!ht]
$$\epsfbox{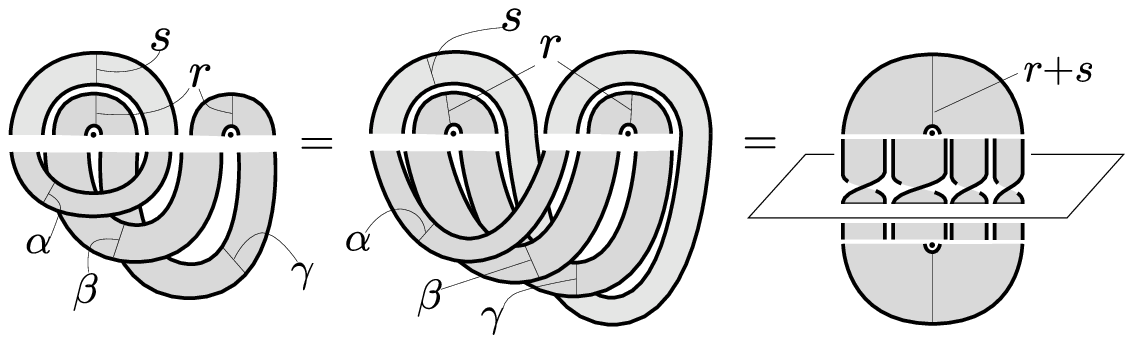}$$
\caption{}\label{fig:equi-4}
\end{figure}

In Figure~\ref{fig:equi-5}, the curve in the middle is obtained
from the one on the left by doubling each components, that is,
taking the boundaries of a two dimensional untwisted regular
neighborhood of the curve on the left. Since there is one arc
component (joining two puctures) in the curve on the left, the
curve in the middle contains $2|(r,s,t,\rho)_{+1}|-1$ components.
Finally the curve in the right is the quotient of the curve in the
middle by the $\mathbf Z/2\mathbf Z$-action given by the
$\pi$-rotation about the axis indicated in the figure. Each
component obtained by doubling a closed component has the linking
number zero with the axis and the component obtained by doubling
the arc component has the linking number 1 with the axis. Thus the
number of components in the curve on the right exactly equals
$|(r,s,t,\rho)_{+1}|$.

\begin{figure}[!ht]
$$\epsfbox{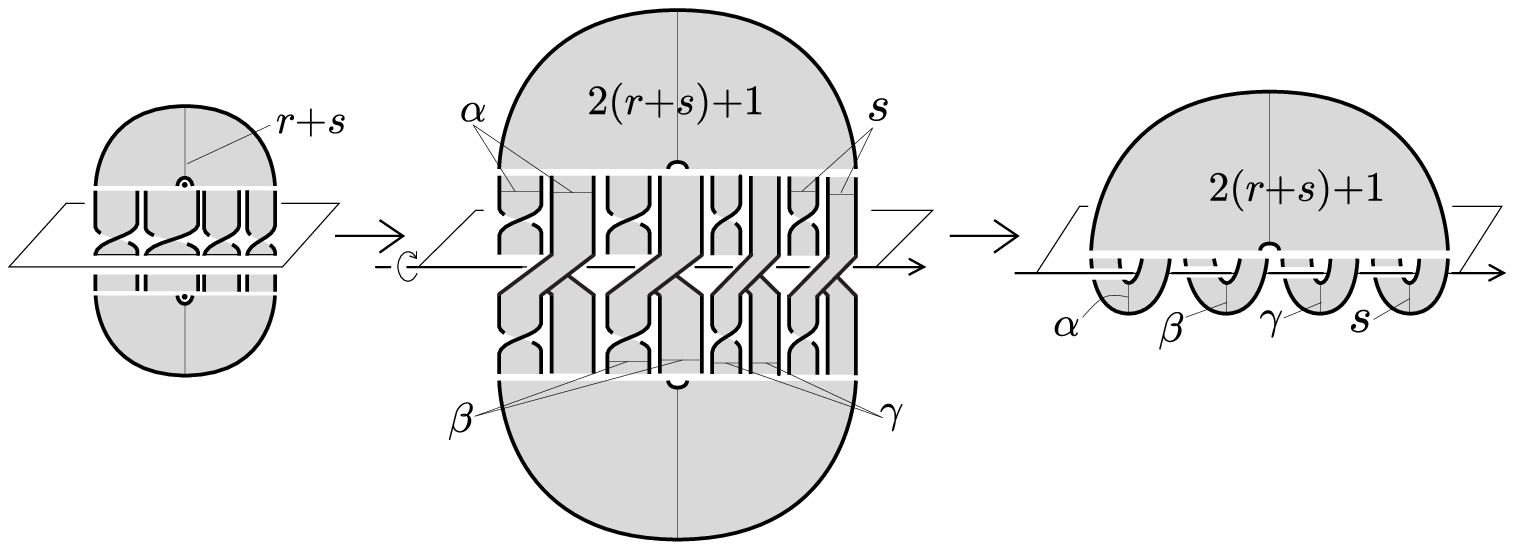}$$
\caption{}\label{fig:equi-5}
\end{figure}
\end{proof}

\begin{lem}\label{lem:3leaves}
In Figure~\ref{fig:23leaves}, the numbers of components in the
curves (a) and (b) are $\gcd(\alpha,\beta)$ and
$\gcd(|\alpha-\gamma|,\beta+\gamma)$, respectively.
\begin{figure}[!ht]
$$\epsfbox{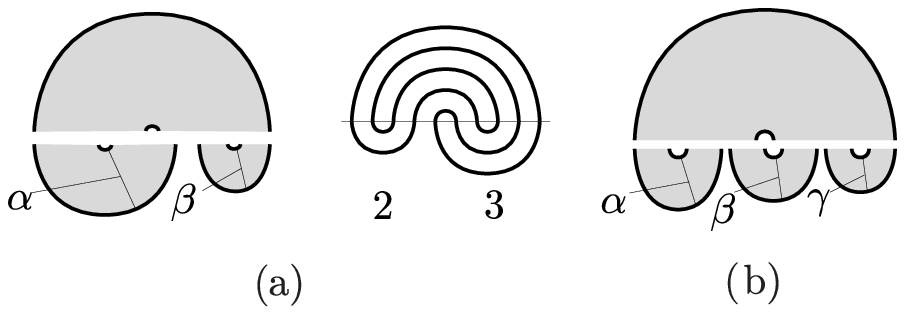}$$
\caption{}\label{fig:23leaves}
\end{figure}
\end{lem}

\begin{proof} It is easy to see that Euclidean algorithm determines
the number of components in the curve (a). The curve (b) can be
modified to a curve like (a) as in the following figure:
$$\epsfbox{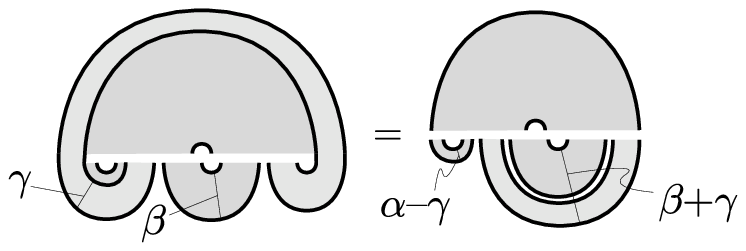}$$
\end{proof}

\begin{figure}[!ht]
$$\epsfbox{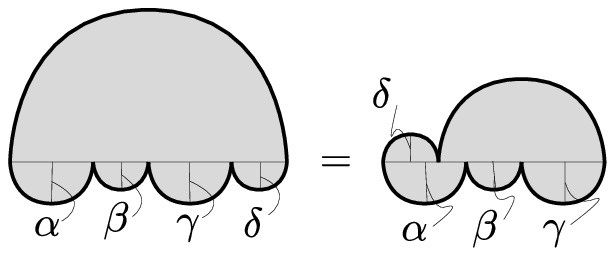}$$
\caption{}\label{fig:4leaves}
\end{figure}

We now consider the curve given in Figure~\ref{fig:4leaves}.

We may assume $\alpha\ge\delta$ by flipping the curve if not. If
$\alpha=\delta$, then it is easy to see that there are $\alpha+
\gcd(\beta,\gamma)$ components in the curve by
Lemma~\ref{lem:3leaves}. If $\alpha>\delta$, then we let
$\bar{\beta}$ be the remainder in the division of $\delta$ by
$\alpha-\delta$, $\bar{\alpha}=\alpha-\delta-\bar{\beta}$,
$\bar{\gamma}=\beta$, and $\bar{\delta}=\gamma$. Via the
modifications in the following figure
$$\epsfbox{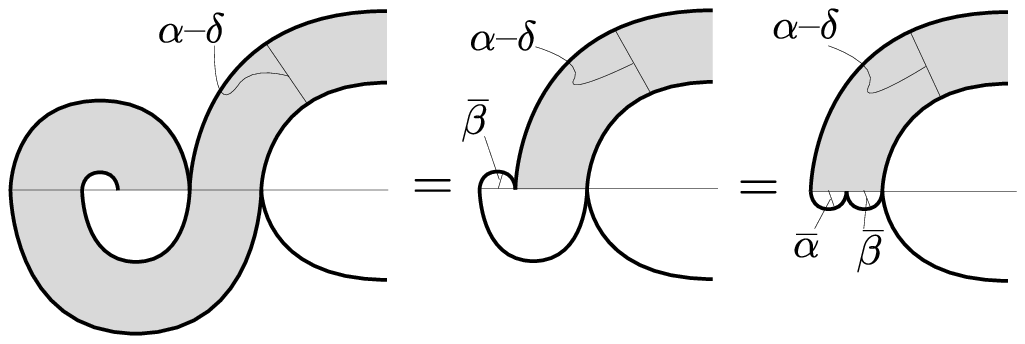}$$
we have
$$\epsfbox{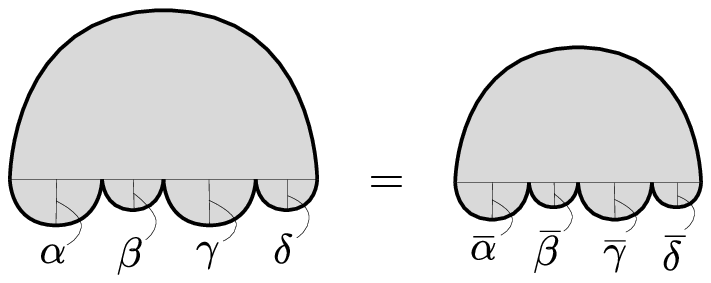}$$
and $\bar\alpha+\bar\beta+\bar\gamma+\bar\delta <
\alpha+\beta+\gamma+\delta$. By repeating this process, we
eventually reach the situation in which at least one of four
numbers becomes zero and so Lemma~\ref{lem:3leaves} can be
applied. Thus we obtain an algorithm to compute
$|(r,s,t,\rho)_{+1}|$ in Proposition~\ref{prop:cc-method} that
equals the number of components of the curve in
Figure~\ref{fig:comp2}. We call this algorithm the \emph{component
counting algorithm}.

We conclude this section by summarizing our result in a theorem.

\begin{thm}\label{thm:class}
Let $(r,s,t,\rho)_\epsilon$ be a 4-tuple of integers such that
$r,s,t\geq 0$ and $\epsilon=\pm 1$. And set $n=2r+1+s+t$ and
$\bar{\rho}\equiv\epsilon\rho\mod n$, $0\leq\bar{\rho}<n$. If
$0\le \rho < n-t$ then we let $\beta\equiv t \mod \rho$ with
$0\le\beta<\rho$, $\alpha=\rho-\beta$, $\gamma=n-\rho-t$, and
$\delta=s$, or if $n-t\le \rho < n$ then we let $\beta\equiv t
\mod (n-\rho)$ with $0\le\beta<n-\rho$, $\alpha=\rho-t$,
$\gamma=n-\rho-\beta$, and $\delta=s$. If the components counting
algorithm returns ``1'' for the input
$(\alpha,\beta,\gamma,\delta)$ then there exists a simple arc with
ends $x,y$ on $T$ which is determined by $(r,s,t,\rho)_\epsilon$.
\end{thm}

\section{Schubert's normal forms and knot groups}\label{sec:normal}

\subsection{Schubert's normal forms of 1-bridge torus knots}
\begin{figure}[!ht]
$$\epsfbox{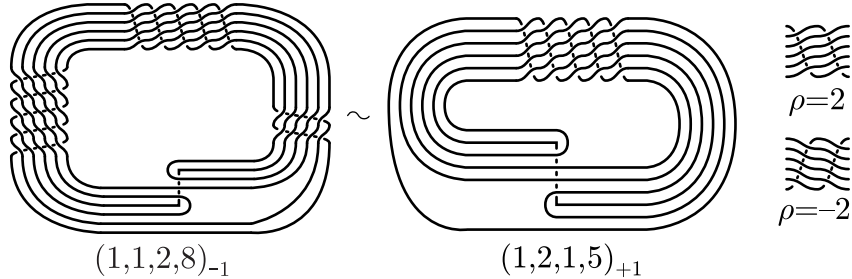}$$
\caption{1-bridge torus knot}\label{fig:1b-ex}
\end{figure}

Let $(r,s,t,\rho)_\epsilon$ be integers satisfying the assumption
in Theorem~\ref{thm:class}. Then  we obtain a simple arc $\alpha$
on $T$ with ends $x,y$. Let $\beta$ be an arc connecting $x$, $y$
in a meridian disk of the solid torus $V_1$ bounded by $T$ and let
$h:\partial V_2\to
\partial V_1$ be a homeomorphism for another solid torus $V_2$.
Then $\alpha\cup\beta$ forms a 1-bridge torus knot in the
3-manifold $M_h$ where $M_h=V_1\cup_h V_2$. We note that a
1-bridge torus knot is embedded in $V_1$ in our convention.
Figure~\ref{fig:1b-ex} contains an example of 1-bridge torus knots
in $S^3$.

Conversely let $(V_1,t_1)\cup_h(V_2,t_2)$ be a (1,1)-decomposition
of a 1-bridge torus knot $K$. If we isotope the trivial arc $t_2$
on the $\partial V_2$ $=$$\partial V_1$ along a cancelling disk
for $t_2$ then it is a simple arc in the torus and so it is
determined by a 4-tuple $(r,s,t,\rho)_\epsilon$ by
Theorem~\ref{thm:arc-class} since the isotopies of type (I) and
(II) do not change the knot type and the given 3-manifold.

Thus a 4-tuple $(r,s,t,\rho)_\epsilon$ together with a
homeomorphism $h$ uniquely determines a 1-bridge torus knot in
$M_h$. This representation will be called a \emph{Schubert's
normal form} and the 1-bridge torus knot in $M_h$ represented by a
Schubert's normal form $(r,s,t,\rho)_\epsilon$ will be denoted by
$S(r,s,t,\rho)_\epsilon$.

The following facts can be observed immediately.
\begin{enumerate}
\item $S(r,s,t,\rho)_{+1}$ is equivalent to $S(r,t,s,\rho+(2r+1))_{-1}$ (See Figure~\ref{fig:1b-ex}).
\item Two 1-bridge torus knots $S(r,s,t,\rho)_{+1}$ and
$S(r,s,t,-\rho)_{-1}$ are  mirror images each other.
\item $S(0,s,t,\rho)_{+1}$ in $S^3$ is a 1-bridge braid studied in \cite{gabai}.
\item A $(p,q)$-torus knot is a 1-bridge torus knot $S(0,0,p-1,-q)_{+1}$ or
$S(0,p-1,0,-q+1)_{-1}$ in $S^3$.
\item Any 2-bridge knot in $S^3$ represented by an ``original'' Schubert's normal form
$B(\alpha,\epsilon\beta)$ is a 1-bridge torus knot
$S(\beta-1,\alpha-2\beta+1,0,\epsilon)_\epsilon$ (See
Figure~\ref{fig:2bvs1bt}), where $\alpha$ and $\beta$ are positive
odd integers with $\beta<\alpha$ and are relatively prime (For
example, see Chapter 3 of \cite{burde}).
\begin{figure}[!ht]
\epsfbox{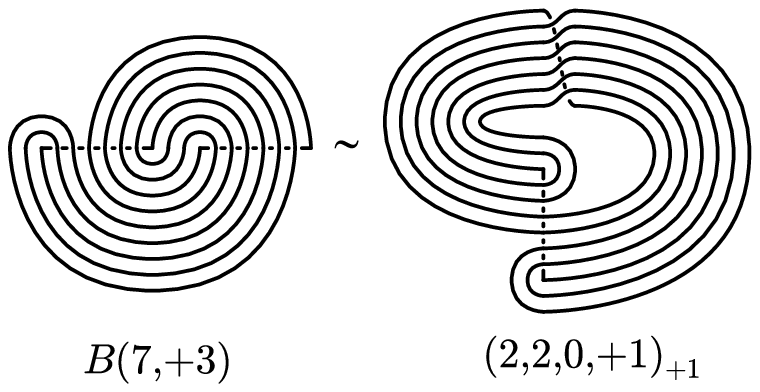}
\caption{}\label{fig:2bvs1bt}
\end{figure}
\item Morimoto and Sakuma showed in
\cite{mori-saku} that any satellite knot which admits an
unknotting tunnel is equivalent to a knot
$K(\alpha,\epsilon\beta;p,q)$ represented by 4 integers such that
$\alpha$ is even and p,q are positive and $0<\beta<\alpha/2$. The
knot $K(\alpha,\epsilon\beta;p,q)$ is a 1-bridge torus knot
$S(\frac{\beta-1}{2},\frac{\alpha-2\beta}{2}, \frac{\alpha}{2}p,
\frac{\alpha}{2}q)_\epsilon.$
\end{enumerate}

In \cite{god-hay}, Goda and Hayashi studied the intersection of
two Heegaard surfaces of $(2,0)$-knots $(M,K)$ that give (1,1) and
(2,0) decomposition, respectively. Using \cite{kob-sae}, they
characterized the intersection when $M$ admits a double cover
branched along $K$. In the next section we will give a necessary
and sufficient condition for the existence of this double cover in
terms of parameters in a Schubert's normal form.

\subsection{Exteriors of 1-bridge torus knots}
In this section we compute the fundamental groups of exteriors of
1-bridge torus knots in a 3-manifold. As a corollary, we give a
necessary and sufficient condition that a lens space has a double
cover branched along a given 1-bridge torus knot in terms of
parameters of its Schubert's normal form.

\begin{figure}[!ht]
$$\epsfbox{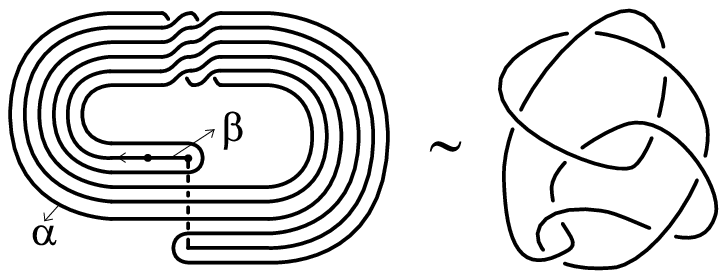}$$
\caption{}\label{fig:9-42}
\end{figure}

\begin{figure}[!ht]
$$\epsfbox{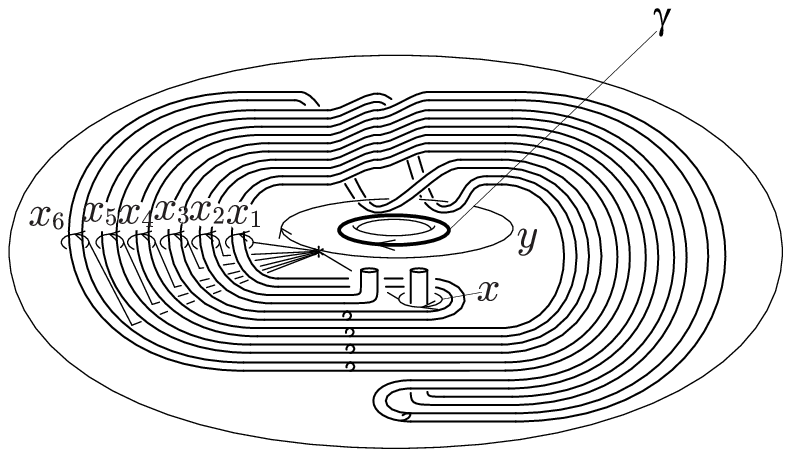}$$
\caption{}\label{fig:ext-9-42}
\end{figure}

Consider a 1-bridge torus knot $K=S(r,s,t,\rho)_{+1}$ in a
3-manifold $M_h$ with a Heegaard decomposition $V_1\cup_h V_2$ of
genus one. We may assume that the homeomorphism $h :
\partial V_2\to \partial V_1$ send the meridian of $\partial V_2$
to the $(p,q)$-curve of $\partial V_1$ and either $(p,q)=(0,1)$ or
$(p,q)=(1,0)$ or $p$ and $q$ are coprime such that $2\leq q<p$.
From now on, the 3-manifold $M_h$ will be denoted by $L(p,q)$.
Thus $L(1,0)=S^3$, $L(0,1)=S^1\times S^2$, and $L(p,q)$ for $1\le
q<p$ is a usual lens space. Let $K$ be expressed as a union of
arcs $\alpha$ and $\beta$ such that $\beta$ is an insignificant
part of $K$ drawn in Figure~\ref{fig:9-42}. We may assume $\alpha$
and $\beta$ are properly embedded in $V_1$ and $V_2$,
respectively. A tubular neighborhood $N(K)$ of $K$ is also the
union $N(\alpha)\cup N(\beta)$ of tubular neighborhoods of
$\alpha$ and $\beta$. Then the exterior $E(K)=L(p,q)-\inte N(K)$
of $K$ is homeomorphic to a 3-manifold obtained from the handle
body $(V_1-\inte N(\alpha))$ by attaching a 2-handle along the
(p,q)-curve determined by $h$(see Figure~\ref{fig:ext-9-42}).

In order to help to understand our algorithm that computes the
fundamental group of $E(K)$, we give an example first.
Let $K$ be the 1-bridge torus knot $S(1,3,0,2)_{+1}$ in $S^3$
shown in Figure~\ref{fig:9-42} which is a mirror image of the knot
$9_{42}$ in the knot table of the Rolfsen's book \cite{rolfsen}.
Since $V_1-\inte N(\alpha)$ is a handle body of genus 2, its
fundamental group is the free group of two generators $x$ and $y$
where $x$, $y$ represent a meridian of $K$ and the core of the
torus as shown in Figure~\ref{fig:ext-9-42} and the base point is
placed on the core. To present the free group $\pi_1(V_1-\inte
N(\alpha))$, we introduce 6 extra generators $x_1,\ldots,x_6$ as
in Figure~\ref{fig:ext-9-42} that are freely homotopic to either
$x$ or $x^{-1}$ and introduce 6 extra relations so that
$$
\begin{array}{r} \pi_1(V_1-\inte N(\alpha))=\langle\:
x,y,x_1,x_2,x_3,x_4,x_5,x_6\mid x_3=W_1^{-1}xW_1,
x_1=W_2^{-1}x_3^{-1}W_2,\\
x_6=W_3^{-1}x_1W_3, x_5=W_4^{-1}x_6W_4, x_4=W_5^{-1}x_5W_5,
x_2=W_6^{-1}x_4^{-1}W_6 \:\rangle
\end{array}
$$
where $W_1=y^{-1}$,
$W_2=x^{-1}$, $W_3=yx$, $W_4=yx$, $W_5=yx$ and
$W_6=yx^{-1}y^{-1}$. These relations are obtained by sliding $x_i$
toward $x$ along the arc $\alpha$ and can be written down without
looking at a knot diagram. We note that $W_i$'s are always words on $x$, $y$.

The attaching curve $\gamma$ for the 2-handle
can be easily described by a word $R$ on $x_1,\ldots,x_6$ and $y$.
The word $R$ can also be generated by an algorithm. In our example
$R=x_6^{-1}x_5^{-1}y$.
After applying Tietze transformations that eliminate these extra generators and relations,
we obtain a presentation of $\pi_1(E(K))$ with two generators $x$, $y$ and
one relation $R$ as follows:
$$\pi_1(E(K))=\langle\:x,y\mid (W_1W_2W_3)^{-1}x(W_1W_2W_3)
(W_1W_2W_3W_4)^{-1}x(W_1W_2W_3W_4)y\:\rangle.$$

Since a 1-bridge torus knot group is always presented by two generators
and one relation, it is easy to compute its Alexander polynomial
via free differential calculus. In our example, the image of derivative
of $R$ with respect to $x$
under the linear extension $\colon\Z(\pi_1(E(K)))\to \Z[t^{\pm
1}]$ of the abelization map is given by $-t^{-1}(1-t-t^2-t^3-t^4+t^5)$ and
the image of derivative of $R$ with respect to $y$ is given by
$-t(1-2t+t^2-2t^3+t^4)$. Thus the Alexander polynomial of $K$ is
their greatest common divisor $1-2t+t^2-2t^3+t^4$.

We now describe an algorithm to compute the fundamental group of
the exterior of a 1-bridge torus knot $K$ given by its Schubert's
normal form $S(r,s,t,\rho)_{+1}$. Let $n=2r+1+s+t$ and
$\rho=mn+\bar\rho$ for $0\le\bar\rho<n$. As we did in the previous
example, we choose small loops $x_1,\ldots,x_n$ that are freely
homotopic to $x$ or $x^{-1}$ so that the product $x_1\cdots x_n$
represents the $(0,1)$-curve  on the torus. The free group
$\pi_1(V_1-\inte N(\alpha))$ generated by $x$, $y$ is presented by
adding $n$ extra generators $x_1,\ldots,x_n$ and $n$ extra
relations $x_{k_i}=W_i^{-1}x_{k_{i-1}}^{\epsilon_i}W_i$ for
$i=1,\ldots,n$ where $x_{k_0}=x$ and $W_{i}$ is a word on $x$ and
$y$. The permutation $\{k_1,\ldots,k_n\}$ of $\{1,\ldots,n\}$ is
obtained by following the arc $\alpha$ starting from the $x$ end.
The integer $k_i$ should be taken modulo $n$ in the following
algorithm that computes the triple $(k_i, W_i, \epsilon_i)$.

\begin{enumerate}
\item[(I)] $(k_1,W_1,\epsilon_1)=(s+r+1-\bar\rho, y^{-1}, +1)$.
\item[(II)] $(k_i, W_i, \epsilon_i)$ is computed inductively as follows:
\begin{enumerate}
\item[(i)] If $\epsilon_1\cdots\epsilon_{i-1}=+1$, then
$$(k_i,W_{i},\epsilon_i)=
\left\{\begin{array}{ll} (2(r+1)-k_{i-1},x,-1) &,1\leq k_{i-1}\leq r\\
(2(r+1)-k_{i-1},x^{-1},-1) &,r+1<k_{i-1}\leq 2r+1\\
(k_{i-1}-(2r+1)-\bar\rho,(yx)^{-1}, +1)&,2r+1<k_{i-1}\leq n-t\\
(k_{i-1}-\bar\rho,y^{-1},+1)&,n-t<k_{i-1}\leq n
\end{array}\right.
$$
\item[(ii)] If $\epsilon_1\cdots\epsilon_{i-1}=-1$, then for
$j\equiv k_{i-1}+\bar\rho \mod n$ and $1\le j\le n$,
$$(k_i,W_{i},\epsilon_i)=
\left\{\begin{array}{ll}
(j+2r+1,yx, +1)&,1\leq j \leq s\\
(2(r+s+1)-j-\bar\rho,yxy^{-1},-1)&,s<j\leq s+r\\
(2(r+s+1)-j-\bar\rho,yx^{-1}y^{-1},-1)&,s+r+1<j\leq n-t\\
(j,y,+1)&,n-t<j\leq n
\end{array}\right.
$$
\end{enumerate}
\end{enumerate}

We now find the word $R(p,q)$ in $\pi_1(V_1-\inte N(\alpha))$ that
represents the $(p,q)$-curve $\gamma$ on $\partial V_1$. Recall
our restriction that either $(p,q)=(0,1)$ or $(p,q)=(1,0)$ or $p$ and
$q$ are coprime with $0<q<p$.
it is clear that $R(0,1)=x_1x_2\cdots x_n$ and
$$
R(1,0)=\left\{\begin{array}{ll} y&,\rho=0\\
\underbrace{R(0,1)^{-1}\cdots R(0,1)^{-1}}_{m\text{
times}}x_n^{-1}x_{n-1}^{-1}
\cdots x_{n-\bar\rho+1}^{-1}y&,\rho>0\\
\underbrace{R(0,1)\cdots R(0,1)}_{-m\text{
times}}x_n^{-1}x_{n-1}^{-1}\cdots x_{n-\bar\rho+1}^{-1}y&,\rho<0
\end{array}\right.
$$
For a given $q\ge 1$, let $R(p,q)=R(p)$ for $p\ge 1$ and $R(p)$
can be computed inductively as follows:
\begin{enumerate}
\item[(I)] $R(1)=R(1,0)$
\item[(II)]
$R(i)=\left\{\begin{array}{ll} R(i-1)R(1,0)&\text{ if }1\leq j\leq
p-q\\
R(i-1)R(0,1)R(1,0)&\text{ otherwise}
\end{array}\right.
$\\
where $j\equiv 1+(i-1)q \mod p$ and $1\le j\le p$.
\end{enumerate}
Therefore a presentation of $\pi_1(E(K))$ is given by
$$
\langle\: x,y,x_1,\ldots,x_n\mid R(p,q),\:
x_{k_i}=W_{i}^{-1}x_{k_{i-1}}^{\epsilon_i}W_{i},
i=1,\ldots,n\:\rangle
$$
Using the substitutions $x_{k_i}= W_i^{-1}\cdots
W_1^{-1}x^{\epsilon_1\cdots\epsilon_i}W_1\cdots W_i$, the word
$R(p,q)$ is written only on $x$ and $y$ and consequently
$\pi_1(E(K))=\langle\: x,y\mid R(p,q)\:\rangle$ by Tietze
transformations.

We consider the first homology of the exterior $E(K)$ of a
1-bridge torus knot $K$ in $L(p,q)$. Notice that it is obviously
infinitely cyclic if the ambient 3-manifold is $L(1,0)=S^3$.
\begin{coro}\label{coro:e-homology}
Let $K$ be a 1-bridge torus knot in $L(p,q)$.
$$
H_1(E(K))\cong\left\{
\begin{array}{ll}
\Z\oplus \Z_{\gcd(p,\ell)}&, p\ge 1\\
\Z\oplus \Z_{\ell}&, p= 0
\end{array}\right.
$$
where
$\ell=\sum\limits_{j=1}^{n}(\prod\limits_{i=1}^{j}\epsilon_i)$.
\end{coro}
\begin{proof}
Let $n_x(p,q)$ (or $n_y(p,q)$) be the exponent sum of $x$ (or $y$,
respectively) in the relation $R(p,q)$ above. Then $n_x(0,1)=\ell$,
$n_y(0,1)=0$, and $n_y(1,0)=1$. Let $n_x(1,0)=m$. Then
$n_x(p,q)=q\ell+pm$ and $n_y(p,q)=p$.
Since $H_1(E(K))$ is the abelianization of
$\pi_1(E(K))$, it has an abelian presentation $\langle\:
x,y\mid (q\ell+pm)x+py\:\rangle$.
Since $p$ and $q$ are coprime,
$H_1(E(K))$ is isomorphic to $\langle\:
x',y'\mid \ell x'+py'\:\rangle\cong\Z\oplus\Z_{\gcd(p,\ell)}$.
If $p=0$, then $q=1$ and so the result follows.
\end{proof}

Using the first homology of the exterior, it is now easy to
determine when a lens space admits a $k$-fold cyclic branched cover
branched along a 1-bridge torus knot $K$.
In the proof of the above corollary, $H_1(E(K))=\langle\:
x,y\mid (q\ell+pm)x+py\:\rangle$ without changing generators
and $x$ represents a meridian of
$K$. A $k$-fold cyclic branched cover exits if and only if
there is an epimorphism $\phi:H_1(E(K))\to \Z_k$ sending $x$ to 1.
This is equivalent to find $\phi(y)\in \Z$ satisfying
$q\ell+pm+p\phi(y)\equiv 0 \mod k$. Such a $\phi(y)$ exists if and
only if $q\ell+pm$ is a multiple of $\gcd(p,k)$. Since the exponent
sum $m$ of $x$ in the word $R(1,0)$ is not hard to compute, the last
condition is readily verified. If $(p,q)=(0,1)$, then
$H_1(E(K))=\langle\:x,y\mid \ell x\:\rangle$ and so we have the condition
that $\ell$ is a multiple of $k$ instead.
When $k=2$, the last condition is simple to describe in term of
parameters in a Schubert's normal form.

\begin{coro}\label{coro:e-2branch}
Let a 1-bridge torus knot $K$ in a lens space $L(p,q)$ have a
Schubert's normal form $S(r,s,t,\rho)_\epsilon$. $L(p,q)$ has
a double branched cover branched along $K$ if and only if either
$p$ is odd or $p$ is even and $s+t$ is odd.
\end{coro}

\begin{proof} Since $k=2$, $\gcd(p,k)=1$ or 2. If $p$ is odd,
$\gcd(p,k)=1$ and so the last condition above always holds. If $p$
is even, $\gcd(p,k)=2$ and the last condition holds iff $\ell$ is
even since $q$ must be odd. Since $K$ intersects the meridian disk
of $V_1$ geometrically $2r+s+t+1$ times, $2r+s+t+1\equiv \ell \mod
2$ and so $s+t$ is odd. This proof includes the case when
$(p,q)=(0,1)$.
\end{proof}

In \cite{cat-mul}, the first homology of the exterior of a
1-bridge torus knot was given in terms of an abstract description
of a 1-bridge torus knot via a mapping class of twice punctured
torus and also a necessary and sufficient condition for the
existence of $k$-fold branched cover of a lens space along a
1-bridge torus knot was given. However it is hard to apply their
result to a 1-bridge torus knot given by a diagram, for example.

\section{Conway's normal forms and Double branched covers} In
this section, we describe double branched covers of a 3-manifolds
$L(p,q)$ branched along 1-bridge torus knots and compute their
first homologies as invariants of 1-bridge torus knots. In order
to do this, we introduce another parameterization of 1-bridge
torus knots that is an analogue of Conway's normal forms of
2-bridge knots.

\subsection{Conway's normal forms of 1-bridge torus knots}\label{subsec:conway}
Let $M(1,0)$ and $M(1,2)$ denote the mapping class groups
of a torus and a two-punctured torus, respectively.
By ignoring punctures, we obtain a homomorphism $j_*\colon M(1,2)\to
M(1,0)$. Let $h$ be a mapping class in $M(1,2)$ that sends
a $(0,1)$-curve of $\Bd V_2$ to a $(p,q)$-curve of $\Bd V_1$.
Then $(V_1,t_1)\cup_h (V_2,t_2)$ is a
(1,1)-decomposition of a 1-bridge torus knot $K$ in a 3-manifold $L(p,q)$.
Let $\bar t_i$ be an arc on
$\partial V_i$ such that $t_i\cup\bar t_i$ bounds a cancelling
disk in $V_i$ for $i=1,2$. Then
a knot diagram of $K$ can be recovered by either
$t_1\cup h(\bar t_2)$ on the side of $V_1$ or
$t_2\cup h^{-1}(\bar t_1)$ on the side of $V_2$.
Let $h_0$ be a mapping class in $M(1,2)$ such that $j_*(h_0)=j_*(h)$
and $h_0$ fixes a small disk containing two punctures.
Then $\bar h=h\circ h_0^{-1}$ is in $\ker j_*$.
For a fixed ambient space $L(p,q)$, a 1-bridge torus knot
is essentially determined by an element of the subgroup $\ker j_*$.
Using results in \cite{birman1,birman2} by Joan Birman,
we will give a presentation of the group $\ker j_*$ in the following
lemma. Let
$\sigma$ be a homeomorphism exchanging two punctures as
illustrated in Figure~\ref{fig:exchange} and $\tau_m$ (resp.
$\tau_\ell$) be a homeomorphism sliding one of punctures along the
meridian (resp. longitude) as illustrated in
Figure~\ref{fig:gen-m(1,2)}.
\begin{figure}[!ht]
$$\epsfbox{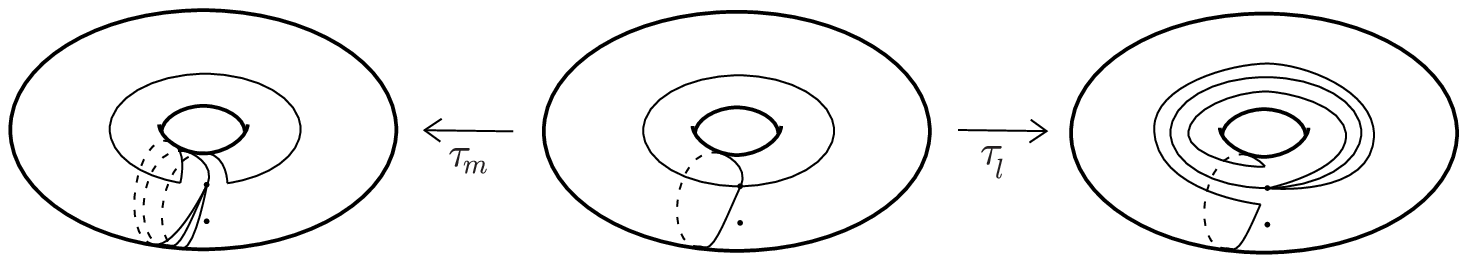}$$
\caption{}\label{fig:gen-m(1,2)}
\end{figure}


\begin{lem}\label{lem:word-rel}
The group $\ker j_*$ is generated by $\tau_\ell$, $\tau_m$ and $\sigma$
under the following defining relations:
\begin{enumerate}
\item $\tau_\ell\sigma^{-1} = \sigma\tau_\ell^{-1}$ and
$\tau_m\sigma = \sigma^{-1}\tau_m^{-1}$,\label{w-rel:1}
\item $\sigma^2=\tau_m^{-1}\tau_\ell^{-1}\tau_m\tau_\ell$.\label{w-rel:2}
\end{enumerate}
\end{lem}
\begin{proof}
The group $\ker j_*$ is generated by $\tau_\ell$, $\tau_m$ and $\sigma$
by Theorem~9 in \cite{birman1} and it is easy to check the above
defining relation from the relations of Theorem 9 in \cite{birman1}
and Corollary 1.3 in \cite{birman2}.
\end{proof}

Since a composition of homeomorphisms (or mapping classes) can be also
regarded as a group multiplication, some
confusion may arises. For example their written orders are opposite
and so homeomorphisms act from the right as an element of a group.
In this paper we will distinguish two operations
 by placing ``$\circ$" between homeomorphisms when they are composed.

\begin{figure}[!ht]
$$\epsfbox{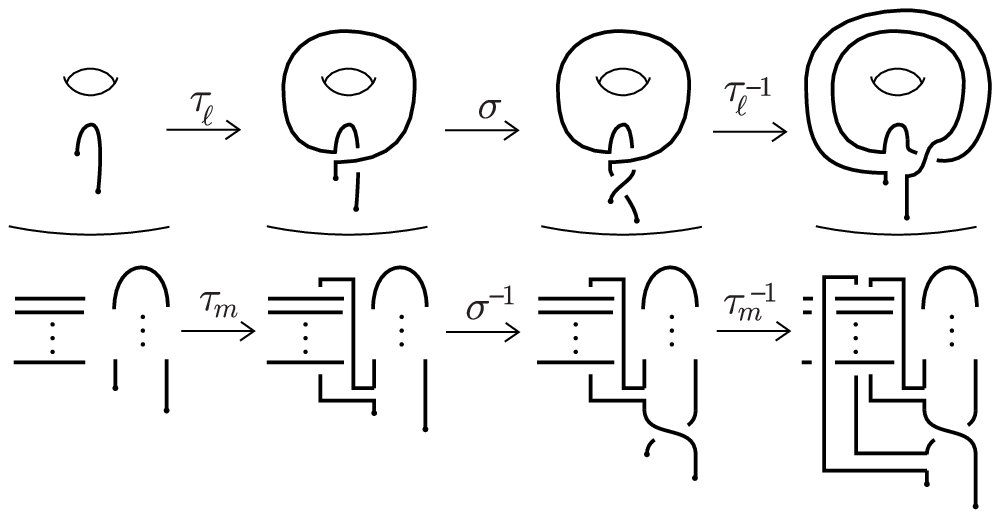}$$
\caption{}\label{fig:effect-hlm}
\end{figure}

Consider the homeomorphisms
$h_\ell=\tau_\ell\sigma^{-1}\tau_\ell^{-1}$ and
$h_m=\tau_m\sigma\tau_m^{-1}$. Figure~\ref{fig:effect-hlm} shows
their effect on a knot by giving the action by $h_\ell^{-1}$
and $h_m^{-1}$ on the arc $t_1$ in $V_1$.

\begin{lem}\label{lem:sub-free}
Let $H$ be the subgroup
of $\ker j_*$ generated by $\sigma$, $h_\ell$, and $h_m$. Then
\begin{enumerate}
\item A word in $\ker j_*$ belongs to $H$ if and only
if both the exponent sum of $\tau_\ell$ and
the exponent sum of $\tau_m$ in the word are even.
\item The subgroup $H$ is a free group.
\end{enumerate}
\end{lem}
\begin{proof}
(1) is clear
since $\sigma\tau_\ell^{-2} = h_\ell$ and $\sigma^{-1}\tau_m^{-2}
= h_m$ by the relation (\ref{w-rel:1}) of Lemma~\ref{lem:word-rel}.

We will show that $H$ is a free group using the
Schreier-Reidemeister method (See \cite{magnus}). By
Lemma~\ref{lem:word-rel}, $X=\{\sigma, \tau_\ell, \tau_m\}$ is a
set of generators and $R=\{(\tau_\ell\sigma^{-1})^2,$
$(\tau_m\sigma)^2,$
$\sigma^2\tau_\ell^{-1}\tau_m^{-1}\tau_\ell\tau_m\}$ is a set of
defining relators for the group $\ker j_*$. The subgroup $H$ is
normal in $\ker j_*$ since conjugates of generators of $H$ by
generator of $\ker j_*$ are again in $H$ as one can see in the
list: $\tau_\ell^{\mathstrut}h_\ell^{\mathstrut}\tau_\ell^{-1}=
h_\ell^{-1}\sigma^{-1}h_\ell^{\mathstrut}$,
$\tau_m^{-1}h_\ell^{\mathstrut}\tau_m^{\mathstrut}=
\sigma^{\mathstrut}h_m^{\mathstrut}h_\ell^{-1}$,
$\tau_m^{\mathstrut}h_\ell^{\mathstrut}\tau_m^{-1}=
h_\ell^{-1}\sigma^{\mathstrut}h_m^{\mathstrut}$,
$\tau_m^{\mathstrut}h_m^{\mathstrut}\tau_m^{-1}=
h_m^{-1}\sigma^{\mathstrut}h_m^{\mathstrut}$,
$\tau_\ell^{-1}h_m^{\mathstrut}\tau_\ell^{\mathstrut}=
\sigma^{-1}h_\ell^{\mathstrut}h_m^{-1}$, and
$\tau_\ell^{\mathstrut}h_m^{\mathstrut}\tau_\ell^{-1}=
h_m^{-1}\sigma^{-1}h_\ell^{\mathstrut}$.
In fact, $\ker j_*/H$ is
isomorphic to $\Z_2\oplus\Z_2$ and so we can choose
$T=\{1,\tau_m, \tau_\ell, \tau_m\tau_\ell \}$ as a Schreier
representative system. By Schreier-Reidemeister method, the
subgroup $H$ is generated by $s$-symbols $s_{k,g}=kg\overline{kg}^{~-1}$ where
$k\in T$, $g\in X$ and $\overline{kg}$ is
a coset representative of $kg$ in $\ker j_*/H$. And $H$ has the defining
relations $t(krk^{-1})$ where $k\in T$, $r\in R$ and $t(\cdot)$ is
a rewriting process. Hence this presentation of $H$ has the following
generators:
\begin{equation*}
\begin{split}
s_{1,\sigma}, s_{1,\tau_\ell}=1, s_{1,\tau_m}=1,s_{\tau_m,\sigma}, s_{\tau_m,\tau_\ell}=1,
s_{\tau_m,\tau_m},\\
s_{\tau_\ell,\sigma}, s_{\tau_\ell,\tau_\ell}, s_{\tau_\ell,\tau_m},
s_{\tau_m\tau_\ell,\sigma}, s_{\tau_m\tau_\ell,\tau_\ell}, s_{\tau_m\tau_\ell,\tau_m}
\end{split}
\end{equation*}
and defining relations:
\begin{align}
\label{relator1}t(\tau_\ell\sigma^{-1}\tau_\ell\sigma^{-1})&
=s_{1,\tau_\ell}^{\mathstrut}s_{\tau_\ell,\sigma}^{-1}s_{\tau_\ell,\tau_\ell}^{\mathstrut}
s_{1,\sigma}^{-1}=s_{\tau_\ell,\sigma}^{-1}s_{\tau_\ell,\tau_\ell}^{\mathstrut}
s_{1,\sigma}^{-1},\\
\label{relator2}t(\tau_m\sigma\tau_m\sigma)&
=s_{1,\tau_m}^{\mathstrut}s_{\tau_m,\sigma}^{\mathstrut}s_{\tau_m,\tau_m}^{\mathstrut}
s_{1,\sigma}^{\mathstrut}=s_{\tau_m,\sigma}^{\mathstrut}s_{\tau_m,\tau_m}^{\mathstrut}
s_{1,\sigma}^{\mathstrut},\\
\label{relator3}t(\sigma^2\tau_\ell^{-1}\tau_m^{-1}\tau_\ell\tau_m)&
=s_{1,\sigma}^2s_{\tau_\ell,\tau_\ell}^{-1}s_{\tau_m\tau_\ell,\tau_m}^{-1}
s_{\tau_m\tau_\ell,\tau_\ell}^{\mathstrut}s_{\tau_m,\tau_m}^{\mathstrut},
\end{align}
\begin{align}
\label{relator4}t(\tau_m\tau_\ell\sigma^{-1}\tau_\ell\sigma^{-1}\tau_m^{-1})&
=s_{1,\tau_m}^{\mathstrut}s_{\tau_m,\tau_\ell}^{\mathstrut}s_{\tau_m\tau_\ell,\sigma}^{-1}
s_{\tau_m\tau_\ell,\tau_\ell}^{\mathstrut}s_{\tau_m,\sigma}^{-1}s_{1,\tau_m}^{-1}
=s_{\tau_m\tau_\ell,\sigma}^{-1}s_{\tau_m\tau_\ell,\tau_\ell}^{\mathstrut}
s_{\tau_m,\sigma}^{-1},\\
\label{relator5}t(\tau_m\tau_m\sigma\tau_m\sigma\tau_m^{-1})&
=s_{1,\tau_m}^{\mathstrut}s_{\tau_m,\tau_m}^{\mathstrut}s_{1,\sigma}^{\mathstrut}
s_{1,\tau_m}^{\mathstrut}s_{\tau_m,\sigma}^{\mathstrut}s_{1,\tau_m}^{-1}
=s_{\tau_m,\tau_m}^{\mathstrut}s_{1,\sigma}^{\mathstrut}s_{\tau_m,\sigma}^{\mathstrut},\\
\label{relator6}t(\tau_m\sigma^2\tau_\ell^{-1}\tau_m^{-1}\tau_\ell\tau_m\tau_m^{-1})&
=s_{1,\tau_m}^{\mathstrut}s_{\tau_m,\sigma}^2s_{\tau_m\tau_\ell,\tau_\ell}^{-1}
s_{\tau_\ell,\tau_m}^{-1}s_{\tau_\ell,\tau_\ell}^{\mathstrut}
=s_{\tau_m,\sigma}^2s_{\tau_m\tau_\ell,\tau_\ell}^{-1}s_{\tau_\ell,\tau_m}^{-1}
s_{\tau_\ell,\tau_\ell}^{\mathstrut},
\end{align}
\begin{align}
\label{relator7}t(\tau_\ell\tau_\ell\sigma^{-1}\tau_\ell\sigma^{-1}\tau_\ell^{-1})&
=s_{1,\tau_\ell}^{\mathstrut}s_{\tau_\ell,\tau_\ell}^{\mathstrut}s_{1,\sigma}^{-1}
s_{1,\tau_\ell}^{\mathstrut}s_{\tau_\ell,\sigma}^{-1}s_{1,\tau_\ell}^{-1}
=s_{\tau_\ell,\tau_\ell}^{\mathstrut}s_{1,\sigma}^{-1}s_{\tau_\ell,\sigma}^{-1},\\
\label{relator8}\begin{split}
t(\tau_\ell\tau_m\sigma\tau_m\sigma\tau_\ell^{-1})&
=s_{1,\tau_\ell}^{\mathstrut}s_{\tau_\ell,\tau_m}^{\mathstrut}s_{\tau_m\tau_\ell,\sigma}^{\mathstrut}
s_{\tau_m\tau_\ell,\tau_m}^{\mathstrut}s_{\tau_\ell,\sigma}^{\mathstrut}s_{1,\tau_\ell}^{-1}\\
&=s_{\tau_\ell,\tau_m}^{\mathstrut}s_{\tau_m\tau_\ell,\sigma}^{\mathstrut}
s_{\tau_m\tau_\ell,\tau_m}^{\mathstrut}s_{\tau_\ell,\sigma}^{\mathstrut},
\end{split}\\
\label{relator9}\begin{split}
t(\tau_\ell\sigma^2\tau_\ell^{-1}\tau_m^{-1}\tau_\ell\tau_m\tau_\ell^{-1})&
=s_{1,\tau_\ell}^{\mathstrut}s_{\tau_\ell,\sigma}^2s_{1,\tau_\ell}^{-1}s_{\tau_m,\tau_m}^{-1}
s_{\tau_m,\tau_\ell}^{\mathstrut}s_{\tau_m\tau_\ell,\tau_m}^{\mathstrut}s_{1,\tau_\ell}^{-1}\\
&=s_{\tau_\ell,\sigma}^2s_{\tau_m,\tau_m}^{-1}s_{\tau_m\tau_\ell,\tau_m}^{\mathstrut},
\end{split}
\end{align}
\begin{align}
\label{relator10}\begin{split}
t(\tau_m\tau_\ell^2\sigma^{-1}\tau_\ell\sigma^{-1}\tau_\ell^{-1}\tau_m^{-1})&
=s_{1,\tau_m}^{\mathstrut}s_{\tau_m,\tau_\ell}^{\mathstrut}
s_{\tau_m\tau_\ell,\tau_\ell}^{\mathstrut}s_{\tau_m,\sigma}^{-1}
s_{\tau_m,\tau_\ell}^{\mathstrut}s_{\tau_m\tau_\ell,\sigma}^{-1}
s_{\tau_m,\tau_\ell}^{-1}s_{1,\tau_m}^{-1}\\
&=s_{\tau_m\tau_\ell,\tau_\ell}^{\mathstrut}s_{\tau_m,\sigma}^{-1}
s_{\tau_m\tau_\ell,\sigma}^{-1},
\end{split}\\
\label{relator11}\begin{split}
t(\tau_m\tau_\ell\tau_m\sigma\tau_m\sigma\tau_\ell^{-1}\tau_m^{-1})&
=s_{1,\tau_m}^{\mathstrut}s_{\tau_m,\tau_\ell}^{\mathstrut}s_{\tau_m\tau_\ell,\tau_m}^{\mathstrut}
s_{\tau_\ell,\sigma}^{\mathstrut}s_{\tau_\ell,\tau_m}^{\mathstrut}
s_{\tau_m\tau_\ell,\sigma}^{\mathstrut}s_{\tau_m,\tau_\ell}^{-1}s_{1,\tau_m}^{-1}\\
&=s_{\tau_m\tau_\ell,\tau_m}^{\mathstrut}s_{\tau_\ell,\sigma}^{\mathstrut}
s_{\tau_\ell,\tau_m}^{\mathstrut}s_{\tau_m\tau_\ell,\sigma}^{\mathstrut}
\end{split}\\
\label{relator12}\begin{split}
t(\tau_m\tau_\ell\sigma^2\tau_\ell^{-1}\tau_m^{-1}\tau_\ell\tau_m\tau_\ell^{-1}\tau_m^{-1})&
=s_{1,\tau_m}^{\mathstrut}s_{\tau_m,\tau_\ell}^{\mathstrut}s_{\tau_m\tau_\ell,\sigma}^{2}
s_{\tau_m,\tau_\ell}^{-1}s_{1,\tau_m}^{-1}s_{1,\tau_\ell}^{\mathstrut}
s_{\tau_\ell,\tau_m}^{\mathstrut}s_{\tau_m,\tau_\ell}^{-1}s_{1,\tau_m}^{-1}\\
&=s_{\tau_m\tau_\ell,\sigma}^{2}s_{\tau_\ell,\tau_m}^{\mathstrut}.
\end{split}
\end{align}
The defining relators (\ref{relator1}) and (\ref{relator7}) induce the same relation,
$s_{\tau_\ell,\tau_\ell}=s_{\tau_\ell,\sigma}s_{1,\sigma}$,
(\ref{relator2}) and (\ref{relator5}) induce
$s_{\tau_m,\tau_m}^{\mathstrut}=s_{\tau_m,\sigma}^{-1}s_{1,\sigma}^{-1}$,
(\ref{relator4}) and (\ref{relator10}) induce
$s_{\tau_m\tau_\ell,\tau_\ell}=s_{\tau_m\tau_\ell,\sigma}s_{\tau_m,\sigma}$,
(\ref{relator12}) induces $s_{\tau_\ell,\tau_m}^{\mathstrut}=s_{\tau_m\tau_\ell,\sigma}^{-2}$,
and finally (\ref{relator8}) and (\ref{relator11}) induce
$s_{\tau_m\tau_\ell,\tau_m}^{\mathstrut}=s_{\tau_m\tau_\ell,\sigma}^{-1}
s_{\tau_\ell,\tau_m}^{-1}s_{\tau_\ell,\sigma}^{-1}
=s_{\tau_m\tau_\ell,\sigma}^{\mathstrut}s_{\tau_\ell,\sigma}^{-1}$.

Thus $s$-symbols $s_{\tau_m,\tau_m}$, $s_{\tau_\ell,\tau_\ell}$, $s_{\tau_\ell,\tau_m}$,
$s_{\tau_m\tau_\ell,\tau_\ell}$ and $s_{\tau_m\tau_\ell,\tau_m}$ together with these
defining relators can be deleted from the presentation for $H$ by Tietze transformations.
Furthermore, the remaining defining
relators are:
\begin{align}
\label{relator3'}(\ref{relator3})&=s_{1,\sigma}^2(s_{1,\sigma}^{-1}s_{\tau_\ell,\sigma}^{-1})
(s_{\tau_\ell,\sigma}^{\mathstrut}s_{\tau_m\tau_\ell,\sigma}^{-1})
(s_{\tau_m\tau_\ell,\sigma}^{\mathstrut}s_{\tau_m,\sigma}^{\mathstrut})
(s_{\tau_m,\sigma}^{-1}s_{1,\sigma}^{-1})=1\\
\label{relator6'}(\ref{relator6})&=s_{\tau_m,\sigma}^2(s_{\tau_m,\sigma}^{-1}
s_{\tau_m\tau_\ell,\sigma}^{-1})(s_{\tau_m\tau_\ell,\sigma}^2)
(s_{\tau_\ell,\sigma}^{\mathstrut}s_{1,\sigma}^{\mathstrut})=s_{\tau_m,\sigma}^{\mathstrut}
s_{\tau_m\tau_\ell,\sigma}s_{\tau_\ell,\sigma}s_{1,\sigma}\\
\label{relator9'}(\ref{relator9})&=s_{\tau_\ell,\sigma}^2(s_{1,\sigma}^{\mathstrut}
s_{\tau_m,\sigma}^{\mathstrut})(s_{\tau_m\tau_\ell,\sigma}^{\mathstrut}s_{\tau_\ell,\sigma}^{-1})
\end{align}
Since the defining relator (\ref{relator6'}) induces the relation
$s_{\tau_m\tau_\ell,\sigma}^{\mathstrut}=s_{\tau_m,\sigma}^{-1}s_{1,\sigma}^{-1}
s_{\tau_\ell,\sigma}^{-1}$, the s-symbol
$s_{\tau_m\tau_\ell,\sigma}$ and the defining relator
(\ref{relator6'}) can be deleted. Finally
$$
(\ref{relator9'})=s_{\tau_\ell,\sigma}^2s_{1,\sigma}^{\mathstrut}s_{\tau_m,\sigma}^{\mathstrut}
(s_{\tau_m,\sigma}^{-1}s_{1,\sigma}^{-1}s_{\tau_\ell,\sigma}^{-1})s_{\tau_\ell,\sigma}^{-1}=1
$$
Consequently no defining relators remain and so $H$ is a free group
generated by three $s$-symbols $s_{1,\sigma}$,
$s_{\tau_m,\sigma}$, and $s_{\tau_\ell,\sigma}$ that correspond to
$\sigma$, $\tau_m^{\mathstrut}\sigma\tau_m^{-1}=h_m$, and
$\tau_\ell^{\mathstrut}\sigma\tau_\ell^{-1}=h_\ell^{-1}$, respectively.
\end{proof}

For integers $a_1,\ldots, a_m, b_1,\ldots,b_m$ and $\delta=0$ or
$1$, let $\bar h$ be a homeomorphism in $H$ defined by
\begin{eqnarray*}
\bar h &=& \tau_\ell^{\delta}(\sigma^{-b_1}h_m^{a_1}\sigma^{-b_2}
h_\ell^{a_2})\cdots(\sigma^{-b_{m-1}}h_m^{a_{m-1}}\sigma^{-b_m}h_\ell^{a_m})\\
&=& (h_\ell^{a_m}\circ\sigma^{-b_m}\circ
h_m^{a_{m-1}}\circ\sigma^{-b_{m-1}})\circ\cdots\circ
(h_\ell^{a_2}\circ\sigma^{-b_2}\circ
h_m^{a_1}\circ\sigma^{-b_1})\circ\tau_\ell^{\delta}.
\end{eqnarray*}
Recall that $h_0$ is a homeomorphism of a torus that fixes a neighborhood
of two punctures and produces a 3-manifold $L(p,q)$.
Then $\bar h$ determines a 1-bridge torus knot $K$ in $L(p,q)$ via a
(1,1)-decomposition $(V_1,t_1)\cup_h(V_2,t_2)$ for $h = \bar
h\circ h_0$. Here the expression
\begin{eqnarray}
[(a_m,b_m,a_{m-1},b_{m-1}),\ldots,(a_4,b_4,a_3,b_3),(a_2,b_2,a_1,b_1)]_\delta.
\label{formula:conway}
\end{eqnarray}
will be called \emph{Conway's normal form} of
the 1-bridge torus knot $K$ in $L(p,q)$. When $\delta=0$,
$\delta$ will be omitted in a Conway's normal form.

\begin{thm}\label{thm:conway}
Every 1-bridge torus knot in a 3-manifold $L(p,q)$ of genus $\le 1$
has a Conway's normal form and it is unique as an element of $H$.
Furthermore, if $p$ is odd then a 1-bridge torus knot in
$L(p,q)$ has a Conway's normal form with $\delta=0$.
\end{thm}
\begin{proof}
An element of $\ker j_*$ belongs to the one of the subsets $H$,
$H\tau_m$, $\tau_\ell H$, and $\tau_\ell H\tau_m$ as we have seen in
the proof of Lemma~\ref{lem:sub-free}. Thus for a homeomorphism $\bar h$
in $\ker j_*$, $\bar h=\tau_\ell^\delta g\tau_m^\epsilon$
for some $g\in H$ and $\epsilon, \delta=0$ or $1$. Two homeomorphisms
$\bar h$ and
$\bar h\tau_m^{-\epsilon}$ determine the same 1-bridge torus
knot and so any 1-bridge torus knot is represented by an element
of $H\cup\tau_\ell H$, that is, it has a Conway's normal form, and
since by Lemma~\ref{lem:sub-free} $H$ is a free group generated by
$\sigma$, $h_\ell$ and $h_m$, a Conway's normal form is unique as an
element of $H$ if it is freely reduced.

Suppose that $p$ is odd. Let $f$ be a homeomorphism in $\ker j_*$ that
slides a puncture once around along a $(p,q)$-curve passing through
the puncture. Then for any element $\bar h$ in $\ker j_*$, $\bar h$
and $\bar h\circ f$ represent the same 1-bridge torus knots in the lens
space $L(p,q)$ since the effect of the action by $f$ on the arc $t_1$
is nullified by the attached 2-handle of $L(p,q)$.
Since the word $f$ as an element of $\ker j_*$ must
contain $p$ $\tau_\ell$'s and $p$ is odd, either $\bar h$ or $f\bar h$
contains an even number of $\tau_\ell$'s. This completes the proof
because any word in $H\cup\tau_\ell H$ that contains an even number
of $\tau_\ell$'s must be in $H$.
\end{proof}

\begin{figure}[!ht]
$$\epsfbox{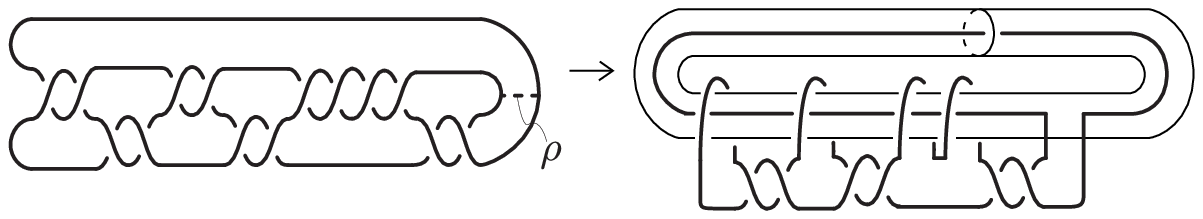}$$
\caption{}\label{fig:conway2b-1bt}
\end{figure}

Consider a 2-bridge knot that has the Conway's normal form
$[2a_1,2a_2,\ldots,2a_m]$ as a 2-bridge knot. Then we obtain its
(1,1)-decomposition using an unknotting tunnel $\rho$ as
illustrated in Figure~\ref{fig:conway2b-1bt} and the homeomorphism
$\bar h$ in $\ker j_*$ of this (1,1)-decomposition can be read as
$(\tau_\ell^{-1}\circ\sigma)\circ
(\sigma^{2a_m}\circ\tau_m^{a_{m-1}})\circ\cdots\circ(\sigma^{2a_2}\circ
\tau_m^{a_1})$. By converting $\bar h$ into an element of $H$
using the relations in Lemma~\ref{lem:word-rel}, we may obtain a
Conway's normal form of a 2-bridge knot as a 1-bridge torus knot.

Now we describe how to obtain a Conway's normal form from a given
Schubert's normal form of a 1-bridge torus knot in $S^3$. Let
$S(r,s,t,mn+\bar\rho)_\epsilon$ be a Schubert's normal form of a
1-bridge torus knot, where $r,s,t$ nonnegative integers,
$-n<\bar\rho<n$, $m$ is an integer, $\epsilon=\pm 1$ and
$n=2r+1+s+t$. If $n$ is even, we may assume that  $m\equiv\bar\rho
\mod 2$ since $mn+\bar\rho=(m+1)n+(\bar\rho-n)$. If $n$ is odd, we
may assume that $\bar\rho$ is even by changing $\epsilon$ since
$S(r,s,t,\rho)_{+1}$ is equivalent to $S(r,t,s,\rho+(2r+1))_{-1}$.
\begin{figure}[!ht]
$$\epsfbox{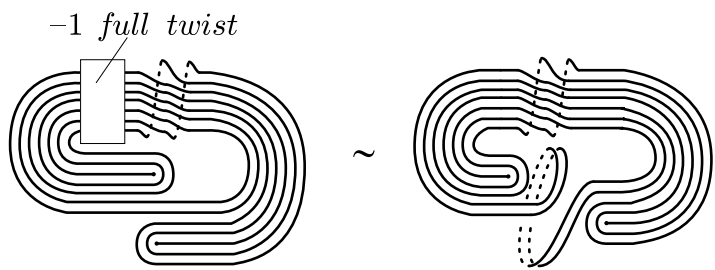}$$
\caption{$S(2,2,0,7-2)_{+1}$}\label{fig:modi-normal}
\end{figure}

\begin{figure}[!ht]
$$\epsfbox{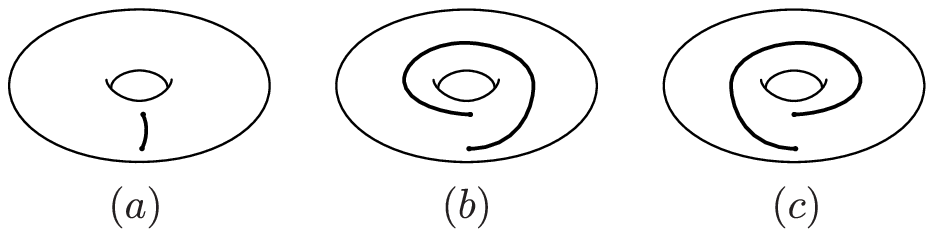}$$
\caption{}\label{fig:2-arcs}
\end{figure}

\begin{figure}[!ht]
$$\epsfbox{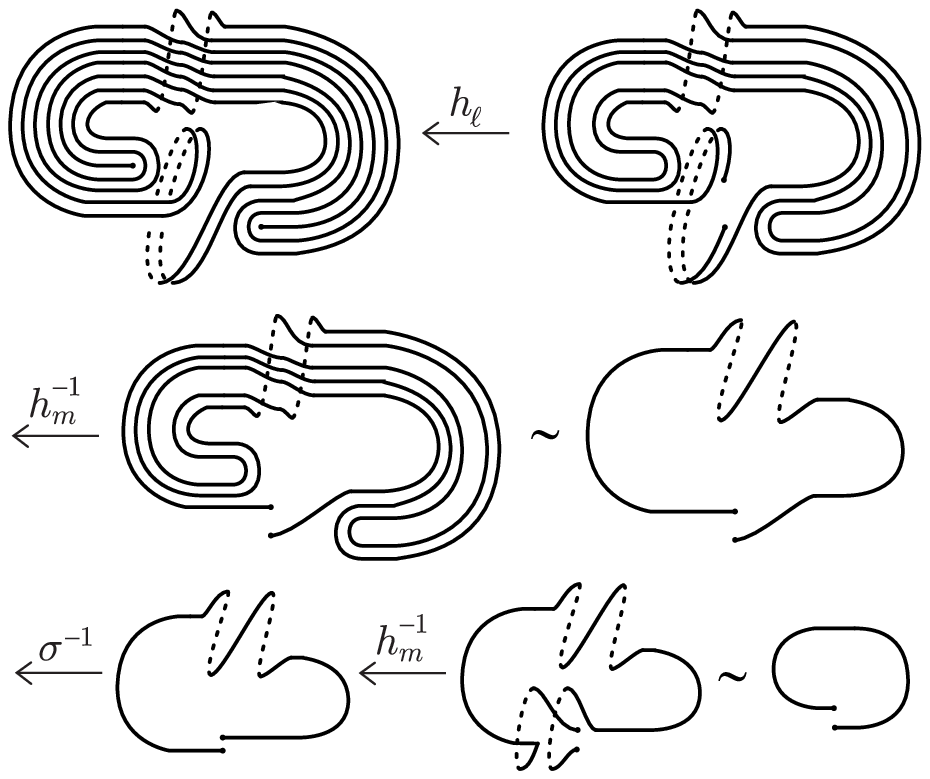}$$
\caption{}\label{fig:normal-conway}
\end{figure}

\begin{figure}[!ht]
$$\epsfbox{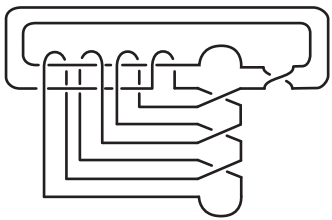}$$
\caption{Conway's normal form $[(1,0,-1,1),(0,0,-1,1)]$ in $S^3$
}\label{fig:conway}
\end{figure}

In order to achieve a symmetry in a knot
$S(r,s,t,mn+\bar\rho)_\epsilon$, we slightly modify the knot
diagram $K$ by moving $m$ full twists in the front as in
Figure~\ref{fig:modi-normal}. Then $K$ intersects the (1,0)-curve
$|m|(s+t)+|\bar\rho|$ times. Under our assumption,
$|m|(s+t)+|\bar\rho|$ is always even. Thus an even number of
$\tau_m$'s must appear in any word $\bar h$ in $\ker j_*$ that
undoes the part $\alpha$ lying on the torus until the number of
intersections between $\alpha$ and the (0,1)-curve becomes 0 or 1
like Figure~\ref{fig:2-arcs}. This means $\bar h$ belongs to the
subgroup $H$. Thus we start to simplify the arc $\alpha$ on the
torus by applying $h_\ell$, $h_m$ and $\sigma$ as in
Figure~\ref{fig:normal-conway} until the diagram becomes one of
three arcs as in Figure~\ref{fig:2-arcs}, all of which form the
trivial knot in $S^3$.
For example, the knot in Figure~\ref{fig:modi-normal} has a
Conway's normal form $[(1,0,-1,1),(0,0,-1,1)]_1$. If the ambient
space is $S^3$ then $[(1,0,-1,1),(0,0,-1,1)]_1$ is isotopic to
$[(1,0,-1,1),(0,0,-1,1)]$ by Theorem~\ref{thm:conway} and it is
isotopic to the knot in Figure~\ref{fig:conway}.

\subsection{Double branched covers along 1-bridge torus knots}
Let $K$ be a 1-bridge torus knot in a 3-manifold $L(p,q)$ with a
Conway's normal form
$$[(a_m,b_m,a_{m-1},b_{m-1}),\ldots,(a_4,b_4,a_3,b_3),
(a_2,b_2,a_1,b_1)]_\delta.$$
We now describe a Heegaard splitting of the double branched cover
of $L(p,q)$ branched along $K$ and compute its first homology.

By Corollary~\ref{coro:e-2branch}, there
exists a double branched cover branched along $K$ if and only if
either $p$ is odd or $p$ is even and the geometric intersection of
$K$ with the meridian disk of $V_1$ is even. Since this geometric
intersection is exactly contributed by the action of each $\tau_\ell$
on $\bar t_2$, the number of $\tau_\ell$'s in $\bar h$ is even and
so $\delta=0$. When $p$ is odd, we may assume $\delta=0$
by Theorem~\ref{thm:conway}. Thus we omit $\delta$ in Conway's normal
forms from now on.

\begin{figure}[!ht]
$$\epsfbox{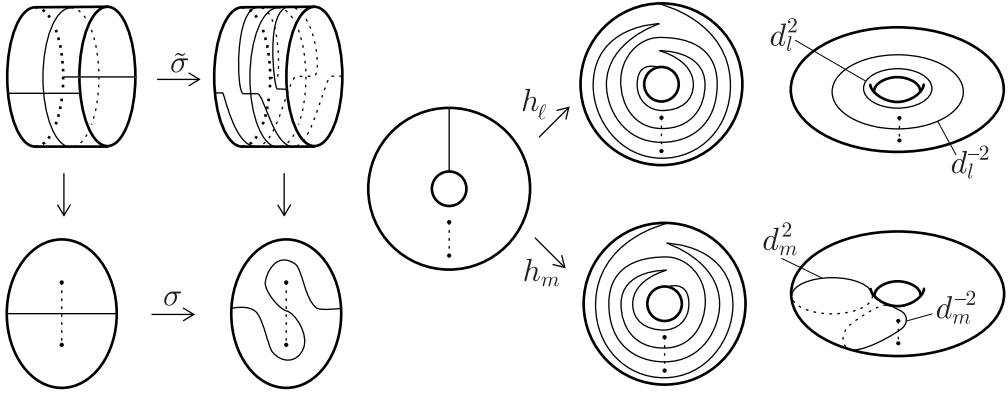}$$
\caption{The lifting of $\sigma$ and the homoemorphisms $h_\ell$,
$h_m$}\label{fig:lift-sh}
\end{figure}

Let a genus two handlebody $\Sigma$ denote the double branched cover
of a solid torus $V$ branched along a trivial arc properly embedded
in $V$.
\begin{figure}[!ht]
$$\epsfbox{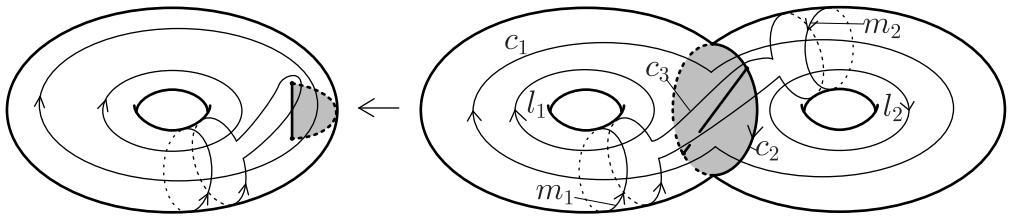}$$
\caption{Double branched cover of a solid torus branched along an arc}
\label{fig:double-st}
\end{figure}
Let $\ell_1,\ell_2,m_1,m_2,c_1,c_2$, and $c_3$ be simple closed
curves on $\partial \Sigma$ as in Figure~\ref{fig:double-st}. Let
$d_c$ denote the Dehn twist along a simple closed curve $c$. Then
the following are evident from Figure~\ref{fig:double-st} and
Figure~\ref{fig:lift-sh}:
\begin{enumerate}
\item The lifting $\tilde h_0$ of $h_0$ is a
homeomorphism of $\Bd\Sigma$ such that $\tilde h_0(m_i)$ (or
$\tilde h_0(\ell_i)$) represents $q[m_i]+p[\ell_i]$ (
$q'[m_i]+p'[\ell_i]$, respectively) on $H_1(\Bd\Sigma)$ for $i=1,2$.
\item The lifting $\tilde\sigma$ of $\sigma$ is $d_{c_2}$.
\item The lifting $\tilde h_\ell$ of $h_\ell$ is
$d_{c_1}^{-1}d_{l_1}^2d_{l_2}^2$.
\item The lifting $\tilde h_m$ of $h_m$ is
$d_{c_3}^{-1}d_{m_1}^2d_{m_1}^2$.
\end{enumerate}

Thus the following theorem is immediate.
\begin{thm}\label{thm:double}
If a 1-bridge torus knot $K$ in a 3-manifold $L(p,q)$ has a Conway's
normal form
$$[(a_m,b_m,a_{m-1},b_{m-1}),\ldots,(a_4,b_4,a_3,b_3),(a_2,b_2,a_1,b_1)]$$
then the double branched cover $X$ of $L(p,q)$ branched along $K$
has a genus-two Heegaard splitting $\Sigma_1\cup_{\tilde
h}\Sigma_2$ where
$$
\tilde h = (\tilde h_\ell^{a_m}\circ\tilde\sigma^{-b_m}\circ\tilde
h_m^{a_{m-1}}\circ\tilde\sigma^{-b_{m-1}})\circ\cdots\circ(\tilde
h_\ell^{a_2}\circ\tilde\sigma^{-b_2}\circ\tilde
h_m^{a_1}\circ\tilde\sigma^{-b_1})\circ\tilde h_0.
$$
\end{thm}

To save notations, $H_1(\Bd\Sigma_1)$ and $H_1(\Bd\Sigma_2)$ are
identified and denoted by $H_1(\Bd\Sigma)$. Then $H_1(\Bd\Sigma)$
is generated by $[m_1]$, $[m_2]$, $[\ell_1]$, and $[\ell_2]$. The
homeomorphism $\tilde\sigma$ induces the identity map on
$H_1(\Bd\Sigma)$. The induced homomorphisms $(\tilde h_\ell^a)_*$
and $(\tilde h_m^b)_*$ are isomorphisms such that
$$
\begin{array}{ll} (\tilde
h_\ell^a)_*([m_1]) = [m_1] + a[\ell_1] -a[\ell_2],& (\tilde
h_\ell^a)_*([\ell_1]) = [\ell_1],\\ (\tilde h_\ell^a)_*([m_2]) =
-a[\ell_1] + [m_2] + a[\ell_2],&(\tilde h_\ell^a)_*([\ell_2]) = [\ell_2],\\
(\tilde h_m^b)_*([m_1]) = [m_1],&(\tilde h_m^b)_*([\ell_1]) =
b[m_1] + [\ell_1] -b[m_2],\\ (\tilde h_m^b)_*([m_2]) =
[m_2],&(\tilde h_m^b)_*([\ell_2]) = -b[m_1] + b[m_2] + [\ell_2].
\end{array}
$$
Therefore $(\tilde h_\ell^a)_*$ and $(\tilde h_m^b)_*$ are
represented by matrices
$$
\left[\begin{array}{cc|cc}
1&0&0&0\\a&1&-a&0\\\hline0&0&1&0\\-a&0&a&1\end{array}\right]
\text{ and } \left[\begin{array}{cc|cc}
1&b&0&-b\\0&1&0&0\\\hline0&-b&1&b\\0&0&0&1\end{array}\right]
\mbox{, respectively.}
$$
Hence the induced homomorphism $\tilde h_*$ is represented by a
matrix $A$ where
$$
A=\left[\begin{array}{cc|cc}1&0&0&0\\a_m&1&-a_m&0\\\hline0&0&1&0\\-a_m&0&a_m&1\end{array}\right]
\cdots\left[\begin{array}{cc|cc}1&0&0&0\\a_2&1&-a_2&0\\\hline0&0&1&0\\-a_2&0&a_2&1\end{array}\right]
\left[\begin{array}{cc|cc}1&a_1&0&-a_1\\0&1&0&0\\\hline0&-a_1&1&a_1\\0&0&0&1\end{array}\right]
\left[\begin{array}{cc|cc}q&q'&0&0\\p&p'&0&0\\\hline0&0&q&q'\\0&0&p&p'\end{array}\right].
$$

\begin{lem}\label{lem:homology}
\begin{align*}
\tilde
h_*([m_1])=&\sum\limits_{i=1}^{2}(-1)^{i-1}\frac{1}{2}(z_m+(-1)^{i-1}q)[m_i]
+(-1)^{i-1}\frac{1}{2}(z_{m+1}+(-1)^{i-1}p)[l_i],
\end{align*}
where $z_m$ is a sequence such that
$z_m=2a_{m-1}z_{m-1}+z_{m-2},~z_0=q\mbox{ and }z_1=p$.
\end{lem}
\begin{proof}
Let
\begin{equation}\label{eqn:tildeh*m1}
\tilde h_*([m_1])=\sum\limits_{i=1}^{2}(x_i[m_i]+y_i[l_i])
\end{equation}
and $A$ be a matrix representing $\tilde h_*$, then
\begin{align*}
A&=\left[\begin{array}{cc|cc}1&0&0&0\\a_m&1&-a_m&0\\\hline0&0&1&0\\-a_m&0&a_m&1\end{array}\right]
\left[\begin{array}{c|c}\begin{array}{c}r_m\\r_{m-1}\\s_m\\s_{m-1}\end{array}&*\end{array}\right]\\
\intertext{and}
\left[\begin{array}{c|c}\begin{array}{c}r_m\\r_{m-1}\\s_m\\s_{m-1}\end{array}&*\end{array}\right]
&=\left[\begin{array}{cc|cc}1&a_{m-1}&0&-a_{m-1}\\0&1&0&0\\\hline0&-a_{m-1}&1&a_{m-1}\\0&0&0&1\end{array}\right]
\left[\begin{array}{c|c}\begin{array}{c}r_{m-2}\\r_{m-1}\\s_{m-2}\\s_{m-1}\end{array}&*\end{array}\right].
\end{align*}
Therefore
\begin{align}
x_1&=r_m,~y_1=a_m(r_m-s_m)+r_{m-1},\label{formula:xy1}\\
x_2&=s_m,~y_2=-a_m(r_m-s_m)+s_{m-1},\label{formula:xy2}\\
r_m&=a_{m-1}(r_{m-1}-s_{m-1})+r_{m-2},\label{formula:rm}\\
s_m&=-a_{m-1}(r_{m-1}-s_{m-1})+s_{m-2},\label{formula:sm}
\end{align}
and $r_0=q$, $r_1=p$ and $s_0=s_1=0$ since
$$
\left[\begin{array}{c|c}\begin{array}{c}r_0\\r_1\\s_0\\s_1\end{array}&*\end{array}\right]=
\left[\begin{array}{c|ccc}q&q'&0&0\\p&p'&0&0\\\hline0&0&q&q'\\0&0&p&p'\end{array}\right].
$$
Let $z_m=r_m-s_m$ and $\bar z_m=r_m+s_m$. Then by the
relations~(\ref{formula:rm})~and~(\ref{formula:sm})
\begin{align*}
z_m&=2a_{m-1}z_{m-1}+z_{m-2},~z_0=q\text{ and }z_1=p.\\
\bar z_m&=\bar z_{m-2},~\bar z_0 = q\text{ and }\bar z_1=p.
\end{align*}
So $\bar z_m=q$ if $m$ is even, $\bar z_m=p$ otherwise. Since $m$ is even,
\begin{align*}
r_m&=1/2(z_m+\bar z_m)=1/2(z_m+q),~r_{m-1}=1/2(z_{m-1}+p),\\
s_m&=-1/2(z_m-\bar z_m)=-1/2(z_m-q),~s_{m-1}=-1/2(z_{m-1}-p).
\end{align*}
Therefore by the
relations~(\ref{formula:xy1})~and~(\ref{formula:xy2})
\begin{align*}
x_1=&r_m=1/2(z_m+q),\\
y_1=&a_mz_m+r_{m-1}=a_mz_m+1/2(z_{m-1}+p)=1/2(z_{m+1}+p),\\
x_2=&s_m=-1/2(z_m-q),\\
y_2=&-a_mz_m+s_{m-1}=-a_mz_m-1/2(z_{m-1}-p)=-1/2(z_{m+1}-p).
\end{align*}
Hence the proof is completed by the
equation~(\ref{eqn:tildeh*m1}).
\end{proof}

For a sequence $\{a_m\}$ and $i=1,2$, define
$$
C_{i,m}^t=\left\{(j_1,\ldots,j_{t})\in \N^t\mid i\leq
j_1<\cdots<j_{t}<m,~j_k-j_{k-1}\geq2,~k=1,\ldots,t\right\},
$$
and
$$
A_i^m(j_1,j_2,\ldots,j_t)= (a_ia_{i+1}\cdots
a_{j_1-1})(a_{j_1+2}\cdots a_{j_2-1})\cdots(a_{j_t+2}\cdots a_m),
$$
Then in particular $A_1^m(1,3,\cdots,m-1)=1$ if $m$ is even and
$A_2^m(2,4,\cdots,m-1)=1$ if $m$ is odd.
Using these notations, we give a solution $z_m$ to the recursive
formula in Lemma~\ref{lem:homology}. But we will omit the proof because
it can be easily observed by an inspection rather than a written proof.
\begin{prop}\label{prop:recursive}
Let $z_m$ be a sequence satisfying the recursive
formula
$$z_{m+1} = 2a_{m}z_{m}+z_{m-1}, z_0=q\mbox{ and }z_1=p,$$
for some sequence $a_m$. Then
\begin{equation}
\label{formula:seq} z_{m+1}= pz_{m+1}^{(1)}+qz_{m+1}^{(2)},
\end{equation}
where for $i=1,2$
\begin{equation}\label{formula:zm+1(1,2)}
z_{m+1}^{(i)}=2^{(m+1)-i}(a_ia_{i+1}\cdots a_m) +
\sum_{t=1}^{[\frac{(m+1)-i}{2}]}2^{(m+1)-i-2t}
\sum_{\substack{(j_1,\ldots,j_{t})\\\in
C^t_{i,m}}}A_i^m(j_1,j_2,\ldots,j_t).
\end{equation}
\end{prop}
From Lemma~\ref{lem:homology} and
Proposition~\ref{prop:recursive}, we calculate the first homology
of $X$.

\begin{thm}\label{thm:dbhomology}
Let $K$ be a 1-bridge torus knot in a 3-manifold $L(p,q)$ with a
Conway's normal form
$$[(a_m,b_m,a_{m-1},b_{m-1}),\ldots,(a_4,b_4,a_3,b_3),(a_2,b_2,a_1,b_1)],$$
and $X$ be a double branched cover of $L(p,q)$ branched along
$K$. Then $H_1(X) \cong \Z_{k_1}\oplus\Z_{k_2}$ where
$k_1=\gcd(p,\frac{1}{2}z_{m+1}^{(2)},z_{m+1})$,
$k_2=\big|pz_{m+1}/k_1\big|$, $z_{m+1}$ is a sequence in
(\ref{formula:seq}), and $z_{m+1}^{(2)}$ is the formula given in
(\ref{formula:zm+1(1,2)}).
\end{thm}
\begin{proof}
$X$ is obtained from $\Sigma_1$ by attaching two 2-handles
along $\tilde h(m_1)$ and $\tilde h(m_1)$ and filling in a 3-ball.
Thus $H_1(X)$ is isomorphic to $H_1(\Sigma_1)$ modulo the image
of $i_*\tilde h_*$ for the inclusion $i:\partial\Sigma_1\to\Sigma_1$.
$H_1(\Sigma_1)$ is a free abelian group
generated by $[\ell_1]$ and $[\ell_2]$.
Lemma~\ref{lem:homology} and the symmetry of $X$,
\begin{align*}
i_*\tilde h_*([m_1])=&\frac{1}{2}(z_{m+1}+p)[\ell_1]-\frac{1}{2}(z_{m+1}-p)[\ell_2],\\
i_*\tilde h_*([m_2])=&\frac{1}{2}(z_{m+1}+p)[\ell_2]-\frac{1}{2}(z_{m+1}-p)[\ell_1].
\end{align*}
Thus $H_1(X)$ is an abelian group with a presentation matrix
$$
R=\left[\begin{array}{rr}
\frac{1}{2}(z_{m+1}+p)&-\frac{1}{2}(z_{m+1}-p)\\
-\frac{1}{2}(z_{m+1}-p)&\frac{1}{2}(z_{m+1}+p)\end{array}\right]\sim
\left[\begin{array}{cc} p&0\\
-\frac{1}{2}(z_{m+1}-p)&z_{m+1}\end{array}\right].
$$
Since $z_{m+1}=pz_{m+1}^{(1)}+qz_{m+1}^{(2)}$ by
Proposition~\ref{prop:recursive},
$$
R\sim
\left[\begin{array}{cc}p&0\\-\frac{1}{2}(p(z_{m+1}^{(1)}-1)+qz_{m+1}^{(2)})&z_{m+1}\end{array}\right]
\sim\left[\begin{array}{cc}p&0\\q\frac{1}{2}z_{m+1}^{(2)}&z_{m+1}\end{array}\right]
$$
Consequently $H_1(X)\cong
\Z_{k_1}\oplus\Z_{k_2}$ where
$$k_1=\gcd(p,q\frac{1}{2}z_{m+1}^{(2)},z_{m+1})=\gcd(p,\frac{1}{2}z_{m+1}^{(2)},z_{m+1})
\text{ and }k_2=\big|pz_{m+1}/k_1\big|.$$
\end{proof}
\begin{coro}\label{coro:dbhomology}
If the ambient space is $S^3$, then $H_1(X)\cong
\Z_{|z_{m+1}^{(1)}|}$ and so
$|z_{m+1}^{(1)}|=|\Delta_K(-1)|$ where $\Delta_K(t)$ is the
Alexander polynomial of $K$.
\end{coro}

\begin{figure}[!ht]
$$\epsfbox{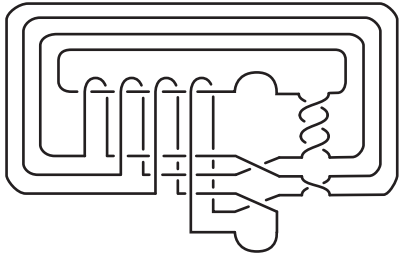}$$
\caption{Conway's normal form
$[(3,0,1,0),(-1,0,1,0)]$}\label{fig:conway-ex}
\end{figure}

The 1-bridge torus knot $K=S(2,2,1,8+1)_{+1}$ in $S^3$ has
a Conway's normal form $[(3,0,1,0),(-1,0,1,0)]$ as in
Figure~\ref{fig:conway-ex}. Therefore
$z_5^{(1)}=2^4(-3)+2^2(-1+3+3)+1 = -27$. By
Corollary~\ref{coro:dbhomology}, we have $H_1(X) = \Z_{27}$. Furthermore
$|\Delta_K(-1)|=27$.

\begin{coro}
Suppose $K$ is a 1-bridge torus knots in $S^3$ with Conway's
normal form
$$[(a_m,b_m,a_{m-1},b_{m-1}),\ldots,(a_4,b_4,a_3,b_3),(a_2,b_2,a_1,b_1)].$$
\begin{enumerate}
\item If $a_{2i} =0$ or $a_{2i-1}=0$ for
$i=1,\ldots,m/2$ then $K$ is a trivial knot.
\item If $a_i$'s are all positive or all negative
then $K$ is a nontrivial knot.
\end{enumerate}
\end{coro}
\begin{proof}
(1) is evident from the knot diagram of $K$. By
Corollary~\ref{coro:dbhomology} the first homology of the double
branched cover is $\Z_{|z_{m+1}^{(1)}|}$. If either $a_i>0$ or
$a_i<0$ for all $i$ then by the formula~(\ref{formula:zm+1(1,2)}),
$z_{m+1}^{(1)}>2^m(a_1a_2\cdots a_m)>1$ since $m$ is even.
Therefore (2) holds.
\end{proof}

\end{document}